\documentclass[12pt,a4paper]{article}
\usepackage{amscd,amsmath,amssymb,amsfonts,times}
\usepackage{mathrsfs}
\usepackage[dvips]{color} 
\usepackage[all]{xy}

\numberwithin{equation}{section}
    \newtheorem{thm}{Theorem}[section]
    \newtheorem{lem}[thm]{Lemma}
    
    \newtheorem{prop}[thm]{Proposition}
    \newtheorem{cor}[thm]{Corollary}

    \newtheorem{rem}[thm]{Remark}
\DeclareMathAlphabet{\mathpzc}{OT1}{pzc}{m}{it}
\newcommand{\qed}
{\mbox{}\nolinebreak$\square$\medbreak\par}
\newenvironment{pf}{\par\smallskip\noindent\emph{Proof.}}{\hfill\qed\par\smallskip}
\newenvironment{pf*}[1]{\par\smallskip\noindent\emph{#1.}}{\hfill\qed\par\smallskip}

\begin{document}
\title{An algorithm of computing special values of 
Dwork's $p$-adic hypergeometric functions in polynomial time}
\author{M. Asakura}
\date\empty
\maketitle

\def\can{\mathrm{can}}
\def\cano{\mathrm{canonical}}
\def\ch{{\mathrm{ch}}}
\def\Coker{\mathrm{Coker}}
\def\crys{\mathrm{crys}}
\def\dlog{{\mathrm{dlog}}}
\def\dR{{\mathrm{d\hspace{-0.2pt}R}}}            
\def\et{{\mathrm{\acute{e}t}}}  
\def\Frac{{\mathrm{Frac}}}
\def\phami{\phantom{-}}
\def\id{{\mathrm{id}}}              
\def\Image{{\mathrm{Im}}}        
\def\Hom{{\mathrm{Hom}}}  
\def\Ext{{\mathrm{Ext}}}
\def\MHS{{\mathrm{MHS}}}  
  
\def\ker{{\mathrm{Ker}}}          
\def\rig{{\mathrm{rig}}}
\def\Pic{{\mathrm{Pic}}}
\def\CH{{\mathrm{CH}}}
\def\NS{{\mathrm{NS}}}
\def\NF{{\mathrm{NF}}}
\def\End{{\mathrm{End}}}
\def\Tr{{\mathrm{Tr}}}
\def\Proj{{\mathrm{Proj}}}
\def\ord{{\mathrm{ord}}}
\def\reg{{\mathrm{reg}}}          %
\def\red{{\mathrm{red}}}          %
\def\Res{\mathrm{Res}}
\def\Spec{{\mathrm{Spec}}}     
\def\syn{{\mathrm{syn}}}
\def\Dw{{\mathrm{Dw}}}
\def\zar{{\mathrm{zar}}}
\def\crys{{\mathrm{crys}}}
\def\logcrys{{\text{\rm log-crys}}}
\def\bA{{\mathbb A}}
\def\bC{{\mathbb C}}
\def\C{{\mathbb C}}
\def\G{{\mathbb G}}
\def\bE{{\mathbb E}}
\def\bF{{\mathbb F}}
\def\F{{\mathbb F}}
\def\bG{{\mathbb G}}
\def\bH{{\mathbb H}}
\def\bJ{{\mathbb J}}
\def\bL{{\mathbb L}}
\def\cL{{\mathscr L}}
\def\bN{{\mathbb N}}
\def\bP{{\mathbb P}}
\def\P{{\mathbb P}}
\def\bQ{{\mathbb Q}}
\def\Q{{\mathbb Q}}
\def\bR{{\mathbb R}}
\def\R{{\mathbb R}}
\def\bZ{{\mathbb Z}}
\def\Z{{\mathbb Z}}
\def\cH{{\mathscr H}}
\def\cD{{\mathscr D}}
\def\cE{{\mathscr E}}
\def\cF{{\mathscr F}}
\def\O{{\mathscr O}}
\def\cV{{\mathscr V}}
\def\cS{{\mathscr S}}
\def\cX{{\mathscr X}}
\def\cY{{\mathscr Y}}
\def\cM{{\mathscr M}}
%
\def\ve{\varepsilon}
\def\vG{\varGamma}
\def\vg{\varGamma}
%
%
%
%
\def\lra{\longrightarrow}
\def\lla{\longleftarrow}
\def\Lra{\Longrightarrow}
\def\hra{\hookrightarrow}
\def\lmt{\longmapsto}
\def\ot{\otimes}
\def\op{\oplus}
\def\l{\lambda}
\def\Isoc{{\mathrm{Isoc}}}
\def\Fil{{\mathrm{Fil}}}
\def\FIsoc{{F\text{-Isoc}}}
\def\FMIC{{F\text{-MIC}}}
\def\Log{{\mathscr{L}{og}}}
\def\FilFMIC{{\mathrm{Fil}\text{-}F\text{-}\mathrm{MIC}}}

\def\wt#1{\widetilde{#1}}
\def\wh#1{\widehat{#1}}
\def\spt{\sptilde}
\def\ol#1{\overline{#1}}
\def\ul#1{\underline{#1}}
\def\us#1#2{\underset{#1}{#2}}
\def\os#1#2{\overset{#1}{#2}}

\def\DF{{{\mathscr DF}}}


\begin{abstract}
Dwork's $p$-adic hypergeometric function is defined to be a ratio
${}_sF_{s-1}(t)/{}_sF_{s-1}(t^p)$ of hypergeometric power series.
Dwork showed that it is a uniform limit of rational functions, and hence one can define
special values on $|t|_p=1$.
However to compute the value modulo $p^n$ in the naive method,
the bit complexity increases by exponential when $n\to\infty$.
In this paper we present a certain algorithm whose complexity increases at most
$O(n^4(\log n)^3)$. 
\end{abstract}

\section{Introduction}
Let $p$ be a prime number.
Let $\ul a=(a_1,\ldots,a_s)\in \Z_p^s$ with $s\geq 2$ an integer.
Let
\[
F_{\ul a}(t):={}_sF_{s-1}\left({a_1,\ldots,a_s\atop 1,\ldots,1};t\right)
=\sum_{n=0}^\infty
\frac{(a_1)_n}{n!}\cdots\frac{(a_s)_n}{n!}t^n\in\Z_p[[t]]
\]
be the hypergeometric series where 
$(\alpha)_n$ denotes the Pochhammer symbol, 
\[
(\alpha)_n:=\alpha(\alpha+1)\cdots(\alpha+n-1),\quad (\alpha)_0:=1.
\]
For $\alpha\in \Z_p$, let $\alpha'$ be the Dwork prime which is defined to
be $(\alpha+k)/p$ with $k\in \{0,1,\ldots,p-1\}$ such that $\alpha+k\equiv0$ mod $p$. 
The ratio
\[
\cF^\Dw_{\ul a}(t):=\frac{F_{\ul a}(t)}{F_{\ul a'}(t^p)},\quad \ul a':=(a'_1,\ldots,a'_s)
\]
is called {\it Dwork's $p$-adic hypergeometric function}.
In his seminal paper \cite{Dwork-p-cycle}, Dwork discovered
a sequence of rational functions which converges $\cF^\Dw_{\ul a}(t)$.
More precisely,
for a power series $f(t)=\sum_{i\geq 0}a_it^i$, we denote by $f(t)_{<m}=[f(t)]_{<m}:=
\sum_{i<m}a_it^i$ the truncated polynomial.
Then Dwork discovered the following congruence relations, 
which we call the Dwork congruences,
\[
\cF^{\Dw}_{\ul a}(t)\equiv 
\frac{F_{\ul a}(t)_{<p^s}}{[F_{\ul a'}(t^p)]_{<p^s}}
\mod p^s\Z_p[[t]].
\]
An immediate consequence of his congruences is that
$\cF^\Dw_{\ul a}(t)$ is a $p$-adic holomorphic function 
in the sense of Krasner (i.e. an element of Tate algebra), 
and hence one can define the special value at $t=\alpha\in \C_p
:=\wh{\ol\Q}_p$
by
\begin{equation}\label{intro-lim}
\cF^{\Dw}_{\ul a}(t)|_{t=\alpha}=\cF^{\Dw}_{\ul a}(\alpha)=
\lim_{n\to \infty}
\left(\frac{F_{\ul a}(t)_{<p^s}}{[F_{\ul a'}(t^p)]_{<p^s}}\bigg|_{t=\alpha}
\right)
\end{equation}
under the condition 
\begin{equation}\label{intro-*}
\bigg|F_{\ul a'}(t)_{<p^n}|_{t=\alpha}\bigg|_p=1,
\quad \forall\,n\geq 1
\end{equation}
where $(-)|_{t=\alpha}$ denote the usual evaluations of rational functions and $|\cdot|_p$ denotes the $p$-adic valuation on $\C_p$. 
Using these values,
he described the Frobenius eigenvalues
on the Monsky-Washnitzer cohomology or {\it rigid cohomology} (cf. \cite{LS}), 
which is now called the unit root formula
(cf. \cite[\S 7]{Put}).
After the work by Dwork, many people studied the congruences,
for example, N. Katz \cite{Katz} developed the congruences in more general situation, 
and moreover alternative methods are brought in
\cite{BV}, \cite{straten}, \cite{MV} etc.

\medskip

On the other hand, 
when we compute a value
\[
\cF^{\Dw}_{\ul a}(t)|_{t=\alpha}\mod p^n
\]
according to the Dwork congruence,
we find a difficulty especially for large $n$.
Indeed the degrees of polynomials increase by
exponential order, and hence the coefficients $A_i,A'_i$ get larger
very quickly,
\[
(a)_{p^n}\sim (p^n)!\sim e^{p^n(n\log p-1)}\quad \text{(Stirling)}.
\]
We note that the bit complexity for computing $(a)_{p^n}$ is $O(n^2p^{2n})$
(by the naive multiplication algorithm).
The aim of this paper is to present a certain algorithm for
computing the special values in case $s=2$, away from the Dwork congruences.

\medskip

\noindent{\bf Main result}.
{\it
Let $N,M\geq 2$ be integers.
Let $a\in \frac{1}{N}\Z$, $b\in \frac{1}{M}\Z$ with $0<a,b<1$. 
Suppose that $p>\max(N,M)$ (hence $p\ne2$).
Let $W=W(\ol\F_p)$ be the Witt ring of $\ol\F_p$. 
Let $\alpha\in W^\times\setminus(1+pW)$ be an arbitrary element
satisfying \eqref{intro-*}.
Then there is an algorithm of computing the special value
\[
\cF^\Dw_{a,b}(t)|_{t=\alpha}\mod p^nW
\]
such that the bit complexity (for fixed $a,b,p,\alpha$) is
at most $O(n^4(\log n)^3)$ as $n\to\infty$.}

\medskip

The algorithm is displayed in \S \ref{algorithm-sect}.

\medskip

Let us see the examples in case that $a=b=1/2$ and 
$p^n=5^{20}\sim9.5367\times 10^{13}$.
It is almost impossible to compute
\[
\left(\frac12\right)_{5^{20}-1},\quad 
\left(\frac12\right)_{5^{20}-2},\ldots
\]
modulo $5^{20}$
by an ordinary PC in a direct way, because they are too large.
On the other hand, our algorithm allows to compute in a few seconds, e.g. 
\begin{align*}
\cF^\Dw_{\frac12,\frac12}(t)|_{t=2}&\equiv7213582472073 \mod 5^{20}\\
\cF^\Dw_{\frac12,\frac12}(t)|_{t=3}&\equiv22359491081212 \mod 5^{20}\\
\cF^\Dw_{\frac12,\frac12}(t)|_{t=4}&\equiv65856465245823\mod 5^{20}\\
&\quad\vdots
\end{align*}

Here is the idea of the main result. 
As is so in \cite{Dwork-p-cycle},
we study the Frobenius endomorphism
on the rigid cohomology rather than the function itself.
The key idea is an application of {\it Kedlaya-Tuitman} \cite{KT}.
It is a fundamental fact in the theory of rigid cohomology that the entries of the Frobenius
matrix are overconvergent.
Kedlaya and Tuitman give effective bounds on overconvergence
(\cite[Theorem 2.1]{KT}).
A noteworthy point is that their bounds increase by polynomial order with respect to $n$,
and this is the key to our main result.

However, our result is not a straightforward application of Kedlaya-Tuitman.
To obtain their bound explicitly (this is necessary in explicit computations),
one has to estimate 
the supremum norm of the entries of the Frobenius matrix.
As long as the author sees, this is a non-trivial and delicate question as
it concerns an integral basis of (log) crystalline cohomology (or
de Rham cohomology with integral coefficients), cf. the proof of Theorem \ref{main-thm}
below.
The author does not know how to get the estimate in a general situation,
while it is enough to consider
the following family in our discussion,
\[
f:X\lra \Spec A=\Spec W[t,(t-t^2)^{-1}],\quad f^{-1}(t)=\{(1-x^N)(1-y^M)=t\}
\]
We call this the {\it hypergeometric fibration}, 
which is indeed a nice family to handle so that
one can solve all the above delicate questions.
The first half of this paper is devoted to a careful computation of the de Rham cohomology
$H^1_\dR(X/A)$ with $A$-coefficient.
Next key step is computations of power series expansions of
the Frobenius matrix where
we follow the method of Lauder \cite{Lauder} (=the deformation method).
However we notice that the power series are centered at the singular fiber 
rather than a smooth fiber,
and then a new technique appears in the computation (e.g. we use the $p$-adic digamma
functions, \cite[\S 2]{New}).
 
\medskip

We hope to obtain a generalization of the algorithm for
$\cF^\Dw_{\ul a}(t)$ with $s\geq 3$, by discussing the rigid cohomology
of a higher dimensional hypergeometric fibration 
\[
(1-x_0^{N_0})\cdots(1-x_d^{N_d})=t,
\]
though I have not worked out.

\medskip

\noindent{\it Acknowledgement}.
I would like to express sincere gratitude to Professor Nobuki Takayama,
to whom I am indebted a lot for the discussion on the bit complexity
of the algorithm.

\section{Dwork's $p$-adic Hypergeometric functions}
Let $p$ be a prime number. 
Let $\Z_p$ be the ring of $p$-adic integers, and $\Q_p$ the fractional field.
Let $\C_p$ be the completion of $\ol\Q_p$. Write $O_{\C_p}=\{|x|_p\leq1\}$
the valuation ring.
\subsection{Definition}\label{DefHG-sect}
For an integer $n\geq 0$, 
we denote by $(\alpha)_n$ the Pochhammer symbol, which is defined by
\[
(\alpha)_n:=\alpha(\alpha+1)\cdots(\alpha+n-1),\quad (\alpha)_0:=1.
\]
Let $s\geq 2$ be an integer.
For $(a_1,\ldots,a_s)\in \Q_p^s$ and $(b_1,\ldots,b_{s-1})\in (\Q_p\setminus
\Z_{\leq0})^{s-1}$, 
the {\it hypergeometric power series} is defined to be
\[
{}_sF_{s-1}\left({a_1,\ldots,a_s\atop b_1,\ldots,b_{s-1}};t\right):=\sum_{n=0}^\infty
\frac{(a_1)_n\cdots(a_s)_n}{(b_1)_n\cdots(b_{s-1})_{n-1}}\frac{t^n}{n!}\in \Q_p[[t]].
\]
In this paper, we only consider the series
\[
F_{\ul a}(t):={}_sF_{s-1}\left({a_1,\ldots,a_s\atop 1,\ldots,1};t\right)
=\sum_{n=0}^\infty
\frac{(a_1)_n}{n!}\cdots\frac{(a_s)_n}{n!}t^n\in\Z_p[[t]]
\]
for $\ul a=(a_1,\ldots,a_s)\in \Z_p^s$.

\medskip

For $\alpha\in \Z_p$, 
let $\alpha'$ denote the Dwork prime, which is defined to be $(\alpha+k)/p$
where $k\in \{0,1,\ldots,p-1\}$ such that $\alpha+k\equiv 0$ mod $p$.
Define the $i$-th Dwork prime by
$\alpha^{(i)}=(\alpha^{(i-1)})'$ and $\alpha^{(0)}:=\alpha$.
Write $\ul a'=(a'_1,\ldots,a'_s)$ and 
$\ul a^{(i)}=(a^{(i)}_1,\ldots,a^{(i)}_s)$.
{\it Dwork's $p$-adic hypergeometric function} is defined to be a power series
\[
\cF^\Dw_{\ul a}(t):=\frac{F_{\ul a}(t)}{F_{\ul a'}(t^p)}\in \Z_p[[t]].
\]
A slight modification is
\[
\cF^{\Dw,\sigma}_{\ul a}(t):=\frac{F_{\ul a}(t)}{F_{\ul a'}(t^\sigma)}\in W[[t]]
\]
for a $p$-th Frobenius $\sigma$ on $W[[t]]$ given by $\sigma(t)=ct^p$, $c\in pW$,
where $W=W(\ol\F_p)$ is the
Witt ring of $\ol\F_p$.

\subsection{Dwork's congruence relations}
In general, neither of the power series $F_{\ul a}(t)$ mod $p\Z_p[[t]]$ or
$\cF^\Dw_{\ul a}(t)$ mod $p\Z_p[[t]]$ terminate
(e.g. \cite[(4.28)]{New}).
Therefore one cannot substitute $t=\alpha\in W$
in $F_{\ul a}(t)$ or $\cF^\Dw_{\ul a}(t)$ directly unless $\alpha\in pW$.
In his seminal paper \cite{Dwork-p-cycle},
Dwork showed that there is a sequence of rational functions which
converges to $\cF^\Dw_{\ul a}(t)$, namely it is a 
{\it $p$-adic analytic function} in the sense of Krasner. 
\begin{thm}[Dwork's congruence relations]\label{D-C-thm}
For a power series $f(t)=\sum_{i\geq 0}a_it^i$, we denote
$f(t)_{<k}=[f(t)]_{<k}=\sum_{0\leq i<k}a_it^i$ the truncated polynomial.
Let $\sigma(t)=ct^p$ with $c\in 1+pW$. 
Then
\begin{equation}\label{D-C-eq0}
\cF^{\Dw,\sigma}_{\ul a}(t)\equiv \frac{F_{\ul a}(t)_{<p^n}}{[F_{\ul a'}(t^\sigma)]_{<p^n}}
\mod p^nW[[t]]
\end{equation}
for any $n\geq 1$.
Hence for $\alpha\in \O_{\C_p}$ satisfying 
\begin{equation}\label{D-C-eq1}
F_{\ul a'}(t)_{<p^n}\big|_{t=\alpha}\not\equiv 0\mod {\frak m}_{\C_p},\quad \forall\,n\geq1
\end{equation}
where ${\frak m}_{\C_p}:=\{|x|_p<1\}$ is the maximal ideal,
one can define a special value of $\cF^\Dw_{\ul a}(t)$ at $t=\alpha$ by
\[
\cF^{\Dw,\sigma}_{\ul a}(t)|_{t=\alpha}=\cF^{\Dw,\sigma}_{\ul a}(\alpha)=\lim_{n\to \infty}
\left(\frac{F_{\ul a}(t)_{<p^n}}{[F_{\ul a'}(ct^p)]_{<p^n}}
\bigg|_{t=\alpha}\right).
\]
\end{thm}
\begin{rem}
One cannot substitute $t=\alpha$ in
$F_{\ul a}(t)$ since it is {\it not} a $p$-adic analytic function.
For example, suppose $\ul a'=\ul a$ and $p\ne2$, the following is wrong !
\[
\cF^{\Dw}_{\ul a}(-1)=\frac{F_{\ul a}(-1)}{F_{\ul a}((-1)^p)}
=\frac{F_{\ul a}(-1)}{F_{\ul a}(-1)}=1.
\]
\end{rem}
\begin{pf}
When $c=1$, this is proven in \cite[p.37, Thm. 2, p.45]{Dwork-p-cycle}.
The general case can be reduced to the case $c=1$ in the following way.
Since
$\cF^{\Dw,\sigma}_{\ul a}(t)=\cF^{\Dw}_{\ul a}(t)\cdot F_{\ul a'}(t^p)/F_{\ul a'}(ct^p)$,
it is enough to show that
\[
\frac{F_{\ul a}(t)}{F_{\ul a}(ct)}\equiv
\frac{F_{\ul a}(t)_{<p^n}}{F_{\ul a}(ct)_{<p^n}}\mod p^{n+1}W[[t]]
\]
in general. Let $F_{\ul a}(t)=\sum_i A_it^i$.
Then the above is equivalent to that
\[
\sum_{i+j=m,i,j\geq0}A_{i+p^n}(c^jA_j)-
A_i(c^{j+p^n}A_{j+p^n})\equiv0\mod p^{n+1}
\]
for any $m\geq 0$. However this is obvious as $c^{p^n}\equiv 1$ mod $p^{n+1}$.
\end{pf}
\begin{cor}
\begin{equation}\label{D-C-eq3}
F_{\ul a}(t)_{<p^n}\equiv F_{\ul a}(t)_{<p}(F_{\ul a'}(t)_{<p})^p
\cdots (F_{\ul a^{(n-1)}}(t)_{<p})^{p^{n-1}}\mod p\Z_p[[t]].
\end{equation}
The condition
\eqref{D-C-eq1} holds if and only if
\begin{equation}\label{D-C-eq2}
F_{\ul a^{(i)}}(t)_{<p}\big|_{t=\alpha}\not\equiv 0\mod {\frak m}_{\C_p},\quad \forall\,i\geq0.
\end{equation}
Moreover we have
\begin{equation}\label{D-C-eq2-1}
\cF_{\ul a}^{\Dw,\sigma}(t)\in 
W[t,h(t)^{-1}]^\wedge:=\varprojlim_{n}\left(W/p^nW[t,h(t)^{-1}]\right),
\quad h(t):=\prod_{i=0}^NF_{\ul a^{(i)}}(t)_{<p}
\end{equation}
with some $N\gg0$.
In particular
this is a $p$-adic analytic function in the sense of Krasner. 
\end{cor}
\begin{pf}
It follows from \eqref{D-C-eq0} that one has
\[
\frac{F_{\ul a}(t)_{<p^n}}{[F_{\ul a'}(t^p)]_{<p^n}}\equiv
\frac{F_{\ul a}(t)_{<p^n}}{(F_{\ul a'}(t)_{<p^{n-1}})^p}\equiv 
F_{\ul a}(t)_{<p}
\mod p\Z_p[[t]].
\]
Then one can show \eqref{D-C-eq3} by induction on $n$.
Notice that a set $\{F_{\ul a^{(i)}}(t)_{<p}\text{ mod }p\}_{i\geq 0}$ 
of polynomials with $\F_p$-coefficients
has a finite cardinal.
Therefore \eqref{D-C-eq2} is a condition for finitely
many $i$'s. 
\eqref{D-C-eq2-1} is now immediate.
\end{pf}
\begin{thm}
Let $f^{(j)}(t)=\frac{d^j}{dt^j}f(t)$ denote the $j$-th derivative. Then
\begin{equation}\label{D-C-eq4}
\frac{F^{(j)}_{\ul a}(t)}{F_{\ul a}(t)}\equiv
\frac{F^{(j)}_{\ul a}(t)_{<p^n}}{F_{\ul a}(t)_{<p^n}}
\mod p^n\Z_p[[t]]
\end{equation}
for all $n\geq 1$.
Hence 
\[
\frac{F^{(j)}_{\ul a}(t)}{F_{\ul a}(t)}
\in W[t,h(t)^{-1}]^\wedge,
\quad h(t):=\prod_{i=0}^NF_{\ul a^{(i)}}(t)_{<p}
\]
is a $p$-adic analytic function in the sense of Krasner, and one can define 
the special values by \eqref{D-C-eq4}.
\end{thm}
\begin{pf}
\cite[p.45, Lem.3.4]{Dwork-p-cycle}.
\end{pf}
\section{Hypergeometric Fibrations}\label{HG-sect}
For a smooth scheme $X$ over a commutative ring $A$, we denote by
$H^*_\dR(X/A):={\mathbb H}^*_\zar(X,\Omega^\bullet_{X/A})$ 
the algebraic de Rham cohomology groups.
\subsection{Setting}
\label{HG-const-sect}
Let $N,M\geq 2$ be an integer.
Let $W$ be a commutative ring such that $NM$ is invertible.
Suppose that $W$ contains a primitive $\operatorname{lcm}(N,M)$-th root of unity.
Later we shall take $W$ to be the Witt ring of a perfect field of characteristic $p$.
Let $\P:=\P^1\times \P^1\times\P^1$
be the product of the projective lines over $W$ with homogeneous coordinates
$(X_0,X_1)\times(Y_0,Y_1)\times (T_0,T_1)$.
We use inhomogeneous coordinates $x:=X_1/X_0$, $y:=Y_1/Y_0$, $t:=T_1/T_0$ 
and $z:=x^{-1}$, $w:=y^{-1}$, $s:=t^{-1}$.
Let $Y_s\subset\P$ be the closed subscheme defined by a homogeneous equation
\[
T_0(X_0^N-X_1^N)(Y_0^M-Y_1^M)=T_1X_0^NY_0^M
\]
over $W$.
Let
\[
\xymatrix{
f_s:Y_s\ar[r]& \P^1=\mathrm{Proj}W[T_0,T_1]
}\]
be the projection onto the 3rd line.
Put $A:=W[t,(t-t^2)^{-1}]$, $U:=\Spec A\subset \P^1$ and 
\[X:=f_s^{-1}(U)=\Spec A[x,y]/((1-x^N)(1-y^M)-t).
\]
Then $X\to U$ is smooth and projective, and a geometric fiber 
is a connected smooth projective
curve of genus $(N-1)(M-1)$ (e.g. the Hurwitz formula).
An open set $Y_s\setminus f_s^{-1}(s=0)$ is smooth over $W$ where 
``$s=0$" denotes the closed subscheme $\Spec W[s]/(s)\subset \P^1$.
There are singular loci $\{s=1-z^N=w=0\}$
and $\{s=z=1-w^M=0\}$ in the affine open set
\[
\Spec W[s,z,w]/(s(1-z^N)(1-w^M)-z^Nw^M)\subset Y_s.
\]
All the singularities are of type ``$xy=z^k$'' where $k=N$ or $k=M$. 
One can resolve them according to Propositions \ref{Appendix-1} in Appendix B.
The fiber $f_s^{-1}(0)$ at $t=0$ is a 
relative simple normal crossing divisor
(abbreviated NCD) over $W$ (see Appendix B for the definition), and all components are $\P^1$.
The fiber $f_s^{-1}(1)$ at $t=1$ is an integral divisor which is smooth
outside the point $(x,y,t)=(0,0,1)$. The normalization of $f_s^{-1}(1)$ is the Fermat
curve $z^N+w^M=1$.
In a neighborhood of the point $(x,y,t)=(0,0,1)$,
the fiber $f_s^{-1}(1)\subset Y_s$ is defined by $x^N+y^M-x^Ny^M=0$ $\Leftrightarrow$
$(x')^N+y^M=0$, $x':=x(1-y^M)^{\frac{1}{N}}$.
One can further resolve it according to Propositions \ref{Appendix-2}
in Appendix B.

Summing up the above, 
we have a smooth projective $W$-scheme $Y$ with a fibration
\[
\xymatrix{
f:Y\ar[r]& \P^1=\mathrm{Proj}W[T_0,T_1]
}\]
which satisfies the following conditions.
Let $D_0:=f^{-1}(0)$,
$D_1:=f^{-1}(1)=\sum_i n_iD_{1,i}$ and $D_\infty:=f^{-1}(\infty)=
\sum_jm_jD_{\infty,j}$ denote the fibers at $\Spec W[t]/(t)$,
$\Spec W[t]/(t-1)$ and $\Spec W[s]/(s)$ respectively.
\begin{enumerate}
\item[(i)]
$f$ is smooth over $U=\Spec W[t,(t-t^2)^{-1}]\subset \P^1$, and $X=f^{-1}_s(U)=
f^{-1}(U)$. 
\item[(ii)]
$D_0$ and $\sum_i D_{1,i}$ and $\sum_jD_{\infty,j}$
are simple 
relative NCD's over $W$.
\item[(iii)]
The multiplicities $n_i$ of $D_1$
are either of $1,iN,jM$ with $i\in \{1,\ldots,M\}$, 
$j\in \{1,\ldots,N\}$.
\item[(iv)]
The multiplicities of $m_j$ of $D_\infty$
are integers $\leq \max(N,M)$.
\item[(v)]
Any components of $D_0$ or $D_\infty$ are $\P^1$.
There is a unique component of $D_1$ which is not $\P^1$. It is
the Fermat curve $z^N+w^M=1$.
\end{enumerate}
Let $\mu_n:=\{\zeta \in W^\times|\zeta^n=1\}$ denote the group of $n$-th roots
of $1$.
For $(\zeta_1,\zeta_2)\in \mu_N\times \mu_M$, the morphism
$(x,y,t)\mapsto (\zeta_1x,\zeta_2y,t)$ extends to an automorphism on $Y$ or $X$, which
we write by $[\zeta_1,\zeta_2]$.

\subsection{$H^1_\dR(X/A)$}\label{XdR-sect}
Let $U_0$ and $U_1$ be the affine open sets of $X$ defined by $X_0Y_0\ne0$
and $X_1Y_1\ne0$ respectively,
\[
U_0=\Spec A[x,y]/((1-x^N)(1-y^M)-t),
\]
\[
U_1=\Spec A[z,w]/((1-z^N)(1-w^M)-tz^Nw^M).
\]
Then $X=U_0\cup U_1$.
For $i\in \{1,\ldots, N-1\}$ and $j\in \{1,\ldots, M-1\}$ let
\begin{equation}\label{H10form}
\omega_{ij}:=N\frac{x^{i-1}y^{j-M}}{1-x^N}dx=-M\frac{x^{i-N}y^{j-1}}{1-y^M}dy
\end{equation}
be rational relative 1-forms on $X/A$.
\begin{lem}\label{H10-lem}
$\omega_{ij}\in \vg(X,\Omega^1_{X/A})$.
\end{lem}
\begin{pf}
Multiplying $x^iy^j$ on 
\[
t(1-t)\frac{y^{-M}}{1-x^N}\frac{dx}{x}=t\frac{dx}{x}-MN^{-1}(1-x^N)\frac{dy}{y}
\]
one sees $\omega_{ij}\in \vg(U_0,\Omega^1_{X/A})$.
Similarly, using an equality
\[
\frac{1}{1-z^N}\frac{dz}{z}=(1-(1-t^{-1})(1-w^M))\frac{dz}{z}-MN^{-1}\frac{dw}{w}
\]
one sees
$\omega_{ij}\in \vg(U_1,\Omega^1_{X/A})$.
\end{pf}
\begin{lem}\label{H01-lem}
Let $H^1(X,\O_X)$ be the Zariski cohomology which is isomorphic to the
cokernel of the Cech complex
\[
\delta:
\vg(U_0,\O_X)\op \vg(U_1,\O_X)\lra \vg(U_0\cap U_1,\O_X),\quad
(u_0,u_1)\longmapsto u_1-u_0.
\]
Write $[f]:=f$ mod $\Image \delta\in H^1(\O_X)$.
Then $H^1(X,\O_X)$ is generated as $A$-module by elements
\[
[x^iy^{j-M}],\quad i\in \{1,\ldots,N-1\},\,j\in \{1,\ldots,M-1\}.
\]
Moreover for any integers $k,l$, there is an element $\alpha\in A$ such that
$[ x^{i+kN}y^{j+lM}]=\alpha[x^{i}y^{j-M}]$ in $H^1(X,\O_X)$.
\end{lem}
\begin{pf}
We first note that if $k,l\leq0$ or $k,l\geq0$ then
$[x^ky^l]=0$ by definition.
Let $i,j$ be integers such that $1\leq i\leq N-1$ and $1\leq j\leq M-1$.
Since $1-t=x^N+y^M-x^Ny^M$, one has
\begin{equation}\label{H01-lem-eq1}
(1-t)^k[x^{i}y^{j-M}]=[x^{kN}\cdot x^iy^{j-M}]=[x^{i+kN}y^{j-M}],\quad \forall\, k\geq0.
\end{equation}
Let $l\geq 1$. Then
$(1-t)x^{i-N}y^{j-lM}=x^{i}y^{j-lM}+x^{i-N}y^{j-(l-1)M}-x^{i}y^{j-(l-1)M}$, and this implies
 \begin{equation}\label{H01-lem-eq2}
[x^iy^{j-lM}]=[x^iy^{j-(l-1)M}],\quad \forall\, l\geq2,
\end{equation}
and for $l=1$
 \begin{equation}\label{H01-lem-eq3}
[x^iy^{j-M}]+[x^{i-N}y^j]=0.
\end{equation}
We claim 
\begin{equation}\label{H01-lem-eq4}
[x^{i+kN}y^{j-lM}]\in A[x^iy^{j-M}],\quad \forall\, k\geq 0,\, l\geq 1.
\end{equation}
If $l=1$, this is nothing other than \eqref{H01-lem-eq1}.
If $k=0$, this follows from \eqref{H01-lem-eq2}.
Suppose $k\geq 1$ and $l\geq 2$. Then
\[
(1-t)[x^{i+(k-1)N}y^{j-lM}]=[x^{i+kN}y^{j-lM}]+[x^{i+(k-1)N}y^{j-(l-1)M}]-[x^{i+kN}y^{j-(l-1)M}].
\]
Hence \eqref{H01-lem-eq4} follows by induction on $k+l$.
In the same way, one can show $[x^{i-kN}y^{j+lM}]\in A[x^{i-N}y^j]$ for all $k\geq 1$ and $l\geq 0$.
Therefore $[x^{i-kN}y^{j+lM}]\in A[x^{i}y^{j-M}]$ by \eqref{H01-lem-eq3}.
This completes the proof.
\end{pf}
\begin{prop}\label{dR-prop1}
\begin{enumerate}
\item[$(1)$]
$ \vg(X,\Omega^1_{X/A})$ is a free $A$-module with basis
\[
\omega_{ij},\quad 
i\in \{1,\ldots,N-1\},\,j\in \{1,\ldots,M-1\}.
\]
\item[$(2)$]
$H^1(X,\O_X)$ is a free $A$-module with basis
\[
[x^{i}y^{j-M}],\quad 
i\in \{1,\ldots,N-1\},\,j\in \{1,\ldots,M-1\}.
\]
\end{enumerate}
\end{prop}
\begin{pf}
For a point $s\in U=\Spec A$, we denote the residue field by $k(s)$, and 
write $X_s:=X\times_A\Spec k(s)$.
Let $q=0,1$.
Since $\dim_{k(s)}H^q(\Omega^{1-q}_{X_s/k(s)})=(N-1)(M-1)$ is constant with respect to $s$,
one can apply \cite[III,12.9]{Ha}, so that 
$H^q(X,\Omega^{1-q}_{X/A})$ is a locally free
$A$-module and the isomorphism $H^q(\Omega^{1-q}_{X/A})\ot k(s)\cong
H^q(\Omega^{1-q}_{X_s/k(s)})$ follows.
Obviously $\omega_{ij}|_{X_s}\ne0$ and they are linearly independent over $k(s)$ since
each $\omega_{ij}$ belongs to the distinct 
simultaneous eigenspace with respect to $\mu_N\times
\mu_M$. Noticing that $\dim H^0(\Omega^1_{X_s/k(s)})=(N-1)(M-1)$, 
one sees that $\{\omega_{ij}|_{X_s}\}_{i,j}$ forms a $k(s)$-basis of $H^0(\Omega^1_{X_s/k(s)})$,
and hence 
that $\{\omega_{ij}\}_{i,j}$ forms a $A$-basis of $H^0(A,\Omega^1_{X/A})$ by Nakayama's lemma. This completes the proof of (1). 
In a similar way, the assertion (2) follows by using Lemma \ref{H01-lem}.
\end{pf}

The algebraic de Rham 
cohomology $H^1_\dR(X/A)$ is described in terms of the Cech complexes.
Let 
\[
\xymatrix{
\vg(U_0,\O)\op\vg(U_1,\O_X)\ar[r]^{d\quad}\ar[d]_\delta&
\vg(U_0,\Omega^1_{X/A})\op
\vg(U_1,\Omega^1_{X/A})\ar[d]^\delta\\
\vg(U_0\cap U_1,\O_X)\ar[r]^d&\vg(U_0\cap U_1,\Omega^1_{X/A})
}
\]
be a commutative diagram where $d$ is the differential map and $\delta$ is given by
$(u_0,u_1)\mapsto u_1-u_0$.
Then the de Rham cohomology $H_\dR(X/A)$ is isomorphic to the cohomology
of the total complex. In particular,
an element of $H^1_\dR(X/A)$ is given as the representative of a cocycle
\[
(f)\times(\omega_0,\omega_1)\in
\vg(U_0\cap U_1,\O_X)\times
\vg(U_0,\Omega^1_{X/A})\op
\vg(U_1,\Omega^1_{X/A})
\]
which satisfies $df=\omega_1-\omega_0$.
Let $\omega_{ij}\in \vg(X,\Omega^1_{X/A})$ be as in Proposition \ref{dR-prop1}.
We denote by the same notation $\omega_{ij}$ the element of $H^1_\dR(X/A)$
via the natural map $\vg(X,\Omega^1_{X/A})\to H^1_\dR(X/A)$, which is the representative of a cocycle
\[
(0)\times (\omega_{ij}|_{U_0},\omega_{ij}|_{U_1}).
\]
We construct a lifting
\begin{equation}\label{H01-cech-eq1}
\eta_{ij}:=(x^iy^{j-M})\times(\eta^0_{ij},\eta^1_{ij})\in H^1_\dR(X/A)
\end{equation}
of $[x^iy^{j-M}]\in H^1(\O_X)$ in the following way.
A dircet computation yields
\begin{equation}\label{H01-cech-eq2}
(j-M)(1-t)x^{i-N}y^{j-M-1}dy
-d(x^iy^{j-M})=
-\left(\frac{(j-M)t}{M}+\frac{i}{N}(1-x^N)\right)\omega_{ij}.
\end{equation}
Note $x^{i-N}y^{j-M-1}dy=-z^{N-i}w^{M-j-1}dw\in \vg(U_1,\Omega^1_{X/A})$,
and the right hand side lies in $\vg(U_0,\Omega^1_{X/A})$ by Lemma \ref{H10-lem}.
Therefore we put
\[
\eta^0_{ij}:=-\left(\frac{(j-M)t}{M}+\frac{i}{N}(1-x^N)\right)\omega_{ij},\quad
\eta^1_{ij}:=-(j-M)(1-t)z^{N-i}w^{M-j-1}dw,
\]
then we get the desired cocycle \eqref{H01-cech-eq1}.
By Proposition \ref{dR-prop1} (2) together with liftings \eqref{H01-cech-eq1}, 
the natural map $H^1_\dR(X/A)\to H^1(\O_X)$
is surjective, and hence one has an exact sequence
\[
0\lra \vg(X,\Omega^1_{X/A})\lra H^1_\dR(X/A)\lra H^1(X,\O_X)\to0.
\]
Thus we get the following theorem.
\begin{thm}\label{XAdR-thm}
$H^1_\dR(X/A)$ is a free $A$-module with basis
\[
\omega_{ij},\,\eta_{ij}\quad 
i\in \{1,\ldots,N-1\},\,j\in \{1,\ldots,M-1\}.
\]
\end{thm}

\subsection{$H^1(\cY,\Omega^\bullet_{\cY/W[[\l]]}(\log D))$}\label{YdR-sect}
Let $\l$ be an indeterminate.
Let $\Spec W[[\l]]\to \P^1$ be the morphism induced by
 $\l=t$, $1-t$ or $t^{-1}$.
Let
\[
\xymatrix{
\cY\ar[r]\ar[d]&Y\ar[d]^f\\
\Spec W[[\l]]\ar[r]&\P^1
}
\]
be the base change. Let $D\subset \cY$ denote the central fiber, namely
$D=D_0,D_1$ or $D_\infty$ by the notation in \S \ref{HG-const-sect}.
The reduced part $D_\red$ is a relative simple NCD over $W$.
Put $\cX:=\cY\setminus D$.
Define a $\O_\cY$-module
\[
\Omega^1_{\cY/W[[\l]]}(\log D):=\Coker\left[\O_\cY\frac{d\l}{\l}\to
\Omega^1_{\cY/W}(\log D)\right]
\]
and consider the cohomology group
\[
H^1(\cY,\Omega^\bullet_{\cY/W[[\l]]}(\log D)):=
H^1_\zar(\cY,\O_\cY\to \Omega^1_{\cY/W[[\l]]}(\log D)).
\]
\begin{prop}\label{YdR-prop1}
If  $N!M!$ is invertible in $W$, 
then $\Omega^1_{\cY/W[[\l]]}(\log D)$ is a locally free $\O_\cY$-module.
\end{prop}
\begin{pf}
If  $N!M!$ is invertible in $W$, 
then each multiplicity of $D$ is invertible in $W$ (see \S \ref{HG-const-sect}).
The assertion can be checked locally on noticing that $f$ is given by
$(x_1,x_2)\mapsto \l=x_1^{r_1}x_2^{r_2}$ with $r_1, r_2$ integers which are 
invertible in $W$.
\end{pf}
\begin{thm}\label{YdR-thm1}
Suppose that $W$ is an integral domain of characteristic zero, and that $N!M!$ is invertible in $W$.
Put
\[
H_\l:=\Image[H^1(\cY,\Omega^\bullet_{\cY/W[[\l]]}(\log D))\to
H^1_\dR(\cX/W((\l)))],
\]
\[
\Fil^1 H_\l:=\Image[\vg(\cY,\Omega^1_{\cY/W[[\l]]}(\log D))\to
H^1_\dR(\cX/W((\l)))].
\]
Then $H_\l$ and $\Fil^1H_\l$ are free $W[[\l]]$-modules of rank $2$ and $1$
respectively. More precisely, the following holds.
\begin{enumerate}
\item[$(1)$]
If $\l=t$, then $\Fil^1H_\l$ has a $W[[\l]]$-basis $\{\omega_{ij}\}$ and
$H_\l$ has a $W[[\l]]$-basis $\{\omega_{ij},\eta_{ij}\}$ where $(i,j)$ runs over the
pairs of
integers such that $1\leq i\leq N-1$ and $1\leq j\leq M-1$.
\item[$(2)$]
If $\l=s=t^{-1}$, then $\Fil^1H_\l$ has a $W[[\l]]$-basis $\{\omega_{ij}\}$ and
$H_\l$ has a $W[[\l]]$-basis $\{\omega_{ij},s\eta_{ij}\}$.
\item[$(3)$]
If $\l=1-t$, set
\[
\omega_{ij}^*:=\begin{cases}
\omega_{ij}&i/N+j/M\geq 1\\
(1-t)\omega_{ij}&i/N+j/M<1,
\end{cases}
\]
\[
\eta_{ij}^*:=\begin{cases}
\eta_{ij}&i/N+j/M\geq 1\\
(1-i/N-j/M)t\omega_{ij}-\eta_{ij}&i/N+j/M<1.
\end{cases}
\]
Then $\Fil^1H_\l$ has a $W[[\l]]$-basis $\{\omega^*_{ij}\}$ and
$H_\l$ has a $W[[\l]]$-basis $\{\omega^*_{ij},\eta^*_{ij}\}$.
\end{enumerate}
\end{thm}
The proof of Theorem \ref{YdR-thm1} shall be given in later sections.
\subsection{Preliminary on Proof of Theorem \ref{YdR-thm1}}\label{YdR-pre-sect}
Let $U_{kl}=\cY\cap\{X_kY_l\ne0\}$, $k,l\in\{0,1\}$ be an affine open set.
Then $\cY=\bigcup_{k=0,1}\bigcup_{l=0,1}U_{kl}$.
The cohomology group $H^i(\cY,\Omega^\bullet_{\cY/W[[\l]]}(\log D))$
is isomorphic to the cohomology of the total complex of the double complex
\[
\xymatrix{
\bigoplus
\vg(U_{ab},\O_\cY)\ar[r]^{d\hspace{1cm}}\ar[d]_\delta&
\bigoplus\vg(U_{ab},\Omega^1_{\cY/W[[\l]]}(\log D))\ar[d]^\delta\\
\bigoplus\vg(U_{ab}\cap U_{cd},\O_\cY)\ar[r]&
\bigoplus\vg(U_{ab}\cap U_{cd},\Omega^1_{\cY/W[[\l]]}(\log D))
}
\]
An element of $H^1(\cY,\Omega^\bullet_{\cY/W[[\l]]}(\log D))$ is represented by
a cocycle
\[
(f_{ab,cd})\times (\alpha_{ab})\in \bigoplus\vg(U_{ab}\cap U_{cd},\O_\cY)\times
\bigoplus\vg(U_{ab},\Omega^1_{\cY/W[[\l]]}(\log D))
\]
which satisfies $f_{ab,ef}=f_{ab,cd}+f_{cd,ef}$ and
\[
\alpha_{cd}|_{U_{ab}\cap U_{cd}}-\alpha_{ab}|_{U_{ab}\cap U_{cd}}=d(f_{ab,cd}).
\]
If we replace $\O_\cY$ with $\O_\cX$ and $\Omega^1_{\cY/W[[\l]]}(\log D)$
with $\Omega^1_{\cX/W((\l))}$ in the above,
we obtain
the algebraic de Rham cohomology group $H^i_\dR(\cX/W((\l)))$.
Let $\omega_{ij},\eta_{ij}\in H^1_\dR(X/A)$ be as in \eqref{H10form} and
\eqref{H01-cech-eq1}.
Then $\omega_{ij}|_{\cX}\in H^1_\dR(\cX/W((\l)))$ is represented by
\[
(0)\times (\omega_{ij}|_{U_{ab}})
\in \bigoplus\vg(U_{ab}\cap U_{cd},\O_\cX)\times
\bigoplus\vg(U_{ab},\Omega^1_{\cX/W((\l))}).
\]
A cocycle
which represents $\eta_{ij}|_\cX$ is given as follows.
We note that
\[
x^{i-N}y^{j-M-1}dy=
-z^{N-i}w^{M-j-1}dw\in \vg(U_{11},\Omega^1_{X/A})
\]
and
\begin{align*}
x^{i-N}y^{j-M-1}dy&=
-x^{i-N}w^{M-j-1}dw\\
&=-(x^N+(1-t)w^M-x^Nw^M)x^{i-N}w^{M-j-1}dw\\
&=-(1-w^M)x^iw^{M-j-1}dw-(1-t)t^{-1}\frac{N}{M}(1-w^M)^2x^{i-1}w^{M-j}dx\\
&\in \vg(U_{01},\Omega^1_{X/A})
\end{align*}
where the 2nd equality follows from $1=x^N+(1-t)w^M-x^Nw^M$.
On the other hand
 $x^{i-N}y^{j-M-1}dy\not \in \vg(U_{00},\Omega^1_{X/A})$
while we have \eqref{H01-cech-eq2}.
Moreover $x^{i-N}y^{j-M-1}dy\not \in \vg(U_{10},\Omega^1_{X/A})$ while we have
\[
(j-M)z^{N-i}y^{j-M-1}dy-(1-t)d(z^{2N-i}y^{j-M})=z^N
\left(\frac{(j-M)t}{M}+\frac{(i-2N)}{N}(1-t)(1-z^N)\right)\omega_{ij}.
\]
Therefore we put
\[
\eta_{ij}^{11}:=-(j-M)(1-t)z^{N-i}w^{M-j-1}dw,\quad
\eta_{ij}^{01}:=-(j-M)(1-t)x^{i-N}w^{M-j-1}dw,
\]
\[
\eta^{00}_{ij}:=-\left(\frac{(j-M)t}{M}+\frac{i}{N}(1-x^N)\right)\omega_{ij},\quad
\]
\begin{align*}
\eta_{ij}^{10}:&=
(1-t)z^N
\left(\frac{(j-M)t}{M}+\frac{(i-2N)}{N}(1-t)(1-z^N)\right)\omega_{ij}\\
&=(1-y^M+z^Ny^M)
\left(\frac{(j-M)t}{M}+\frac{(i-2N)}{N}(1-t)(1-z^N)\right)\omega_{ij}
\end{align*}
and
\[
f_{00,11}=f_{00,01}:=x^iy^{j-M},\quad f_{10,11}=f_{10,01}:=(1-t)^2z^{2N-i}y^{j-M}
=(1-t)(1-y^M+z^Ny^M)z^{N-i}y^{j-M},
\]
\[
f_{01,11}:=0,\quad
f_{00,10}:=x^iy^{j-M}-(1-t)^2z^{2N-i}y^{j-M}=(1-x^N)(x^Ny^M-2x^N-y^M)z^{2N-i}y^j
\]
and $f_{11,00}:=-f_{00,11}$ etc. Then we get a cocycle
\begin{equation}\label{H01-cech-eq3}
(f_{ab,cd})\times (\eta^{ab}_{ij})
\in \bigoplus\vg(U_{ab}\cap U_{cd},\O_\cX)\times
\bigoplus\vg(U_{ab},\Omega^1_{\cX/W((\l))})
\end{equation}
which represents $\eta_{ij}|_\cX\in H^1_\dR(\cX/W((\l)))$.

\subsection{Deligne's canonical extension}\label{de-sect} 
Let $j:\Spec\C((\l))\to \Spec \C[[\l]]$.
Let $(\cH,\nabla)$ be an integrable connection on $\Spec\C((\l))$.
There is a unique subsheaf $\cH_e\subset\cH$ which satisfies the following
conditions
(cf. \cite{zucker} (17)).
\begin{enumerate}
\item[(D1)] $\cH_e$ is a free $\C[[\l]]$-module such that $j^{-1}\cH_e=\cH$,
\item[(D2)]
the connection extends to have log pole, 
$\nabla:\cH_e\to\frac{d\l}{\l}\ot\cH_e$,
\item[(D3)]
each eigenvalue $\alpha$ of $\Res(\nabla)$ satisfies $0\leq \mathrm{Re}(\alpha)<1$, 
where $\Res(\nabla)$ is the map defined by a commutative diagram
\[
\xymatrix{
\cH_e\ar[r]^\nabla\ar[d]&
\frac{d\l}\l\ot\cH_e\ar[d]^{\Res\ot1}\\
\cH_e/\l\cH_e\ar[r]^{\Res(\nabla)}&\cH_e/\l\cH_e.
}
\]
\end{enumerate}
The extended bundle $(\cH_e,\nabla)$ is called {\it Deligne's canonical extension}.

Let $g:V\to\Spec \C[[\l]]$ be a projective flat morphism which is smooth over 
$\Spec \C((\l))$. Let $D$ be the central fiber. Suppose that $D_\red$ is a NCD.
We define a locally free $\O_V$-module
\[
\Omega^1_{V/\C[[\l]]}(\log D):=\Coker\left[\O_V\frac{d\l}{\l}\to
\Omega^1_{V/\C}(\log D)
\right]
\]
and $\Omega^k_{V/\C[[\l]]}(\log D):=\wedge^k\Omega^1_{V/\C[[\l]]}(\log D)$.

Let $U:V\setminus D$ and let 
$(\cH,\nabla)=(H^i_\dR(U/\C((\l))),\nabla)$ be the Gauss-Manin connection on 
$\Spec \C((\l))$.
Then Deligne's canonical extension of $\cH$ is given as follows
(\cite{steenbrink}, (2.18)--(2.20)),
\[
\cH_e\cong H^i(V,\Omega^\bullet_{V/\C[[\l]]}(\log D)).
\]
Moreover $\exp(-2\pi i\Res_P(\nabla))$ agrees with the monodromy
operator on $H_\C=\ker(\nabla^{\mathit an})$ around $\l=0$
(cf. \cite{steenbrink}, (2.21)).

\medskip

We turn to our family $\cY\to \Spec W[[\l]]$.
Let $K:=\Frac(W)$ be the fractional field.
The characteristic of $K$ is zero
by the assumption
in Theorem \ref{YdR-thm1}.
Put $\cY_K:=\cY\times_{W[[\l]]}K[[\l]]$, $\cX_K:=\cX\times_{W[[\l]]}K[[\l]]$
and $D_K:=D\times_WK$.
Let $(H^1_\dR(\cX_K/K((\l))),\nabla)$ be the Gauss-Manin connection
on $\Spec K((\l))$.

\begin{prop}\label{fermat-GM}
Let $\nabla:H^1_\dR(X/A)\to Adt\ot H^1_\dR(X/A)$ be the Gauss-Manin connection.
Then
\[
\begin{pmatrix}
\nabla(\omega_{ij})&\nabla(\eta_{ij})
\end{pmatrix}
=dt\ot\begin{pmatrix}
\omega_{ij}&\eta_{ij}
\end{pmatrix}
\begin{pmatrix}
0&(1-i/N)(1-j/M)\\
(t-t^2)^{-1}&(1-i/N-j/M)(1-t)^{-1}
\end{pmatrix}.
\]
\end{prop}
\begin{pf}
\cite[Proposition 4.15]{New}.
\end{pf}
\begin{prop}\label{Deligne ext-prop1}
Put Deligne's canonical extension
\begin{equation}\label{Deligne-ext}
H_{\l,K}:=H^1(\cY_K,\Omega^\bullet_{\cY_K/K[[\l]]}(\log D_K))
\subset H^1_\dR(\cX_K/K((\l))).
\end{equation}
Then the $K[[\l]]$-basis is given as follows.
\begin{enumerate}
\item[$(1)$]
If $\l=t$, then 
\[
H_{\l,K}=\bigoplus_{i,j}K[[\l]]\omega_{ij}\op K[[\l]]\eta_{ij}
\]
where $(i,j)$ runs over the
pairs of
integers such that $1\leq i\leq N-1$ and $1\leq j\leq M-1$.
\item[$(2)$]
If $\l=s=t^{-1}$, then 
\[
H_{\l,K}=\bigoplus_{i,j}K[[\l]]\omega_{ij}\op K[[\l]]s\eta_{ij}.
\]
\item[$(3)$]
If $\l=1-t$, then 
\[
H_{\l,K}=\bigoplus_{i,j}K[[\l]]\omega^*_{ij}\op K[[\l]]\eta^*_{ij}
\]
where $\omega_{ij}^*$ and $\eta_{ij}^*$ are as in Theorem \ref{YdR-thm1} (3).
\end{enumerate}
\end{prop}
\begin{pf}
The condition (D1) is obvious by Theorem \ref{XAdR-thm}.
It is straightforward from Proposition \ref{fermat-GM}
that (D2) and (D3) are satisfied in each case.
\end{pf}
\subsection{Proof of Theorem \ref{YdR-thm1} (1), (2)}\label{YdR-12-sect}
We prove Theorem \ref{YdR-thm1} in case $\l=t$ and in case $\l=s=t^{-1}$. 
Write $\eta_{ij}^\prime=\eta_{ij}$ in case $\l=t$ and $\eta_{ij}^\prime=\l\eta_{ij}$ in case $\l=s$.
Recall from Theorem \ref{XAdR-thm} the fact that 
\[
H^1_\dR(\cX/W((\l)))\cong W((\l))\ot_AH^1_\dR(X/A)
\]
is a free $W((\l))$-module with basis $\{\omega_{ij},\eta_{ij}; 1\leq i\leq N-1,\,
1\leq j\leq M-1\}$.
It follows from Proposition \ref{Deligne ext-prop1} that
\[
H_\l\subset \bigoplus_{i,j}K[[\l]]\omega_{ij}+K[[\l]]\eta^\prime_{ij}\subset H^1_\dR(\cX_K/K((\l))).
\]
Therefore
\begin{equation}\label{YdR-proof-1}
H_\l\subset \bigoplus_{i,j}W[[\l]]\omega_{ij}+W[[\l]]\eta^\prime_{ij}
\end{equation}
as $K[[\l]]\cap W((\l))=W[[\l]]$.
We show the opposite inclusion, namely
\begin{equation}\label{YdR-proof-2}
\omega_{ij},\,\eta^\prime_{ij}\in H_\l.
\end{equation}
We first show $\omega_{ij}\in H_\l$.
There is an inetger $m\geq 0$ such that
$\l^m\omega_{ij}\in \vg(\cY,\Omega^1_{\cY/W[[\l]]}(\log D))$.
On the other hand $\omega_{ij}\in  \vg(\cY_K,\Omega^1_{\cY_K/K[[\l]]}(\log D))$.
Note that $\Omega^1_{\cY/W[[\l]]}(\log D)$ is a locally free $\O_\cY$-module
(Proposition \ref{YdR-prop1}).
Moreover one can check that the map $a$ in the following diagram is injective. 
\[
\xymatrix{
0\ar[r]&\Omega^1_{\cY/W[[\l]]}(\log D)\ar[r]^{\l^m}\ar[d]
&\Omega^1_{\cY/W[[\l]]}(\log D)\ar[r]\ar[d]&
\Omega^1_{\cY/W[[\l]]}(\log D)/\l^m\ar[r]\ar[d]^a&0\\
0\ar[r]&\Omega^1_{\cY_K/K[[\l]]}(\log D)\ar[r]^{\l^m}&\Omega^1_{\cY_K/K[[\l]]}(\log D)\ar[r]&
\Omega^1_{\cY_K/K[[\l]]}(\log D)/\l^m\ar[r]&0.
}
\]
Therefore we have $\omega_{ij}\in\vg(\cY,\Omega^1_{\cY/W[[\l]]}(\log D))$
by diagram chase.

Next we show $\eta^\prime_{ij}\in H_\l$.
Recall from \eqref{H01-cech-eq3} the cocycle which represents $\eta_{ij}$,
\[
(f_{ab,cd})\times (\eta^{ab}_{ij})
\in \bigoplus\vg(U_{ab}\cap U_{cd},\O_\cX)\times
\bigoplus\vg(U_{ab},\Omega^1_{\cX/W((\l))}).
\]
Therefore it is enough to show
\[
f_{ab,cd}\in \vg(U_{ab}\cap U_{cd},\O_\cY),\quad
\eta^{ab}_{ij}\in \vg(U_{ab},\Omega^1_{\cY/W[[\l]]}(\log D)).
\]
in case $\l=t$, and
\[
s f_{ab,cd}\in \vg(U_{ab}\cap U_{cd},\O_\cY),\quad
s \eta^{ab}_{ij}\in \vg(U_{ab},\Omega^1_{\cY/W[[\l]]}(\log D)).
\]
in case $\l=s$.
However we have shown that $\omega_{ij}\in \vg(\cY,\Omega^1_{\cY/W[[\l]]}(\log D))$.
Thus this is immediate from the explicit descriptions in \S \ref{YdR-pre-sect}.
This completes the proof of \eqref{YdR-proof-2}
and hence Theorem \ref{YdR-thm1} (1), (2).
\subsection{Proof of Theorem \ref{YdR-thm1} (3)}\label{YdR-3-sect}
Let $\l=1-t$. By the same discussion as in \S \ref{YdR-12-sect}, one can show
\begin{equation}
H_\l\subset \bigoplus_{i,j}W[[\l]]\omega^*_{ij}+W[[\l]]\eta^*_{ij},
\end{equation}
and hence it is enough to show 
\begin{equation}\label{YdR-proof-3}
\omega^*_{ij},\,\eta^*_{ij}\in H_\l.
\end{equation}
If $i/N+j/M\geq 1$, then the same discussion as the proof of \eqref{YdR-proof-2}
works.
Suppose that $i/N+j/M<1$.
The same discussion still works for showing $\omega^*_{ij}=\l\omega_{ij}\in H_\l$.
The rest is to show that
\begin{equation}\label{YdR-proof-4}
\eta^*_{ij}=(1-i/N-j/M)t\omega_{ij}-\eta_{ij}
\in H_\l.
\end{equation}
Recall from \eqref{H01-cech-eq3} the cocycle which represents $\eta_{ij}$,
\[
(f_{ab,cd})\times (\eta^{ab}_{ij})
\in \bigoplus\vg(U_{ab}\cap U_{cd},\O_\cX)\times
\bigoplus\vg(U_{ab},\Omega^1_{\cX/W((\l))}).
\]
Hence
\[
(f_{ab,cd})\times ((\eta^*_{ij})^{ab}):=
(f_{ab,cd})\times ((1-i/N-j/M)t\omega_{ij}-\eta^{ab}_{ij})
\]
represents $\eta_{ij}^*\in H^1_\dR(\cX/W((\l)))$.
Each $f_{ab,cd}$ obviously belongs to $\vg(U_{ab,cd},\O_\cY)$.
Therefore it is enough to show that each $(\eta^*_{ij})^{ab}\in
\vg(U_{ab},\Omega^1_{\cY/W[[\l]]}(\log D))$.

\begin{align*}
(\eta^*_{ij})^{11}&=\left(1-\frac{i}{N}-\frac{j}{M}\right)t\omega_{ij}+(j-M)(1-t)z^{N-i}w^{M-j-1}dw,\\
(\eta_{ij}^*)^{01}&=\left(1-\frac{i}{N}-\frac{j}{M}\right)t\omega_{ij}+(j-M)(1-t)x^{i-N}w^{M-j-1}dw,\\
(\eta_{ij}^*)^{10}&=\left(1-\frac{i}{N}-\frac{j}{M}\right)t\omega_{ij}
-(1-t)z^N
\left(\frac{(j-M)t}{M}+\frac{(i-2N)}{N}(1-t)(1-z^N)\right)\omega_{ij},\\
(\eta^*_{ij})^{00}&=
\left(\frac{(j-M)t}{M}+\frac{i}{N}(1-x^N)\right)\omega_{ij}+
\left(1-\frac{i}{N}-\frac{j}{M}\right)t\omega_{ij}
-\frac{i}{N}x^N\omega_{ij}\\
&=\frac{i}{N}(1-t-2x^N)\omega_{ij}.
\end{align*}
Multiplying $Nz^{N-i}y^j$ on an equality
\[
t\frac{y^{-M}}{1-z^N}\frac{dz}{z}=\frac{M}{N}(1-t)(1-z^N)\frac{dy}{y}+t\frac{dz}{z}
\]
one has 
\[
\omega_{ij}=M(1-t)(1-z^N)z^{N-i}y^{j-1}dy+Ntz^{N-i-1}y^jdz.
\]
The shows $\omega_{ij}\in \vg(U_{10},\Omega^1_{\cY/W[[\l]]}(\log D))$.
Similarly, using equalities
\[
t\frac{x^{-N}}{w^M-1}\frac{dw}{w}=-\frac{N}{M}(1-t)(1-w^N)\frac{dx}{x}-t\frac{dw}{w}
\]
\[
\frac{1}{1-z^N}\frac{dz}{z}=(1-(1-t^{-1})(1-w^M))\frac{dz}{z}-\frac{M}{N}\frac{dw}{w},
\]
one has 
\[
\omega_{ij}\in \vg(U_{10}\cup U_{01}\cup U_{11},\Omega^1_{\cY/W[[\l]]}(\log D)).
\]
We thus have $(\eta^*_{ij})^{ab}\in 
\vg(U_{ab},\Omega^1_{\cY/W[[\l]]}(\log D))$ for
$(a,b)=(0,1),(1,0)$ and $(1,1)$.
The rest is the case $(a,b)=(0,0)$, namely we show
\[
\omega^*_{ij}=(1-t)\omega_{ij},\,x^N\omega_{ij}\in \vg(U_{00},\Omega^1_{\cY/W[[\l]]}(\log D)),
\]
(note that $\omega_{ij}$ no longer belongs to 
$\vg(U_{00},\Omega^1_{\cY/W[[\l]]}(\log D))$).
However the former is already shown in \eqref{YdR-proof-3}, and
the latter follows from an equality
\[
x^N\omega_{ij}=-Mx^iy^{j-1}\frac{dy}{1-y^M}=-Mt^{-1}(1-x^N)x^iy^{j-1}dy.
\]
This completes the proof of \eqref{YdR-proof-4}
and hence Theorem \ref{YdR-thm1} (3).

\section{Rigid cohomology and Dwork's $p$-adic Hypergeometric functions}
Let $W=W(k)$ be the Witt ring of a perfect field $k$ of characteristic $p>0$.
Let $A$ be a faithfully flat $W$-algebra.
We mean by a $p^n$-th Frobenius on $A$ an endomorphism $\sigma$ 
such that $\sigma(x)\equiv x^{p^n}$ mod $pA$ for all $x\in A$ and that 
$\sigma$ is compatible with the 
$p^n$-th Frobenius on $W$. We also write $x^\sigma$ instead of $\sigma(x)$.

For a $W$-algebra $A$ of finite type,
we denote by $A^{\dag}$ the weak completion (cf. \cite[p.5]{LP}).
Namely if $A=W[T_1,\cdots,T_n]$, then 
$A^\dag=W[T_1,\cdots,T_n]^\dag$
is the ring of power series $\sum a_\alpha T^\alpha$ such that
for some $r>1$,
$|a_\alpha|r^{|\alpha|}\to0$ as $|\alpha|\to\infty$, and 
if $A=W[T_1,\cdots,T_n]/I$, then
$A^\dag=W[T_1,\cdots,T_n]^\dag/IW[T_1,\cdots,T_n]^\dag$.

\subsection{Rigid cohomology}
Let $W=W(k)$ be the Witt ring of a perfect field $k$ of characteristic $p>0$.
Put $K:=\Frac(W)$ the fractional field.
For a flat $W$-scheme $V$, we denote $V_K:=V\times_WK$ and
$V_k:=V\times_Wk$. For a flat $W$-ring $A$, we denote $A_K:=A\ot_WK$ and
$A_k:=A\ot_Wk$ as well.

\medskip

Let $A$ be a smooth $W$-algebra, and $X$ a smooth $A$-scheme. 
Thanks to the theory due to Berthelot et al, the {\it rigid cohomology groups}
\[
H^*_\rig(X_k/A_k)
\]
are defined.
We refer the book \cite{LS} for the general theory of rigid cohomology.
Here we list the required properties.
Let $A^\dag$ be the weak completion of $A$, and
$A^\dag_K:=A^\dag\ot_WK$. We fix a $p$-th Frobenius $\sigma$
on $A^\dag$. 
\begin{itemize}
\item
$H^*_\rig(X_k/A_k)$ is a finitely generated 
$A^\dag_K$-module.
\item(Frobenius)
The $p$-th Frobenius $\Phi$ on $H^*_\rig(X_k/A_k)$ (depending on $\sigma$)
is defined in a natural way. This is a $\sigma$-linear endomorphism :
\[
\Phi(f(t)x)=\sigma(f(t))\Phi(x),\quad \mbox{for }x\in H^*_\rig(X_k/A_k),\, 
f(t)\in A_K^\dag.
\]
\item(Comparison with de Rham cohomology)
There is the comparison isomorphism with the algebraic de Rham cohomology,
\[
H^*_\rig(X_k/A_k)\cong H^*_\dR(X_K/A_K)
\ot_{A_K} A^\dag_K.
\]
\item(Comparison with crystalline cohomology)
Let $\alpha$ be a $W$-rational point of $\Spec A$ (i.e. a $W$-homomorphism
$\alpha:A\to W$).
Let $X_\alpha:=X\times_{A,\alpha}W$ denote the fiber at $\alpha$.
There is the comparison isomorphism with the crystalline cohomology,
\[
H^*_\rig(X_k/A_k)\ot_{A^\dag_K,\alpha}K\cong 
H^*_\crys(X_{\alpha,k}/W)\ot\Q.
\]
If $\alpha$ satisfies $\sigma^{-1}({\frak m}_\alpha A^\dag_K)=
{\frak m}_\alpha A^\dag_K$ where ${\frak m}_\alpha\subset A$ denotes the ideal
defining $\alpha$, then $\Phi_\alpha:=\Phi$ mod ${\frak m}_\alpha A^\dag_K$ agrees
with the $p$-th Frobenius on the crystalline cohomology.
\end{itemize}
Let $\cY$ be a proper flat scheme over $W[[t]]$ which is smooth over $W((t))$.
Let the central fiber $D$ at $t=0$. Put $\cX:=\cY\setminus D$.
Suppose that $D_\red$ is a relative NCD over $W$
and the multiplicities of components of $D$ are prime to $p$.
Then there is the comparison isomorphism with the log crystalline cohomology 
with log pole $D$ (\cite[Theorem (6.4)]{Kato-log-crys}),
\begin{equation}\label{crys-dR-eq1}
H^*_\logcrys((\cY_{\ol\F_p},D_{\ol\F_p})/(W[[t]],(t)))
\cong H^*_\zar(\cY,\Omega^\bullet_{\cY/W[[t]]}(\log D)).
\end{equation}
Fix a $p$-th Frobenius $\wh\sigma$ on $W[[t]]$ given by $\wh\sigma(t)=ct^p$ with some 
$c\in  1+pW$. 
Then the $\wh\sigma$-linear $p$-th Frobenius $\Phi_\crys$ on the crystalline
cohomology group is defined in a natural way.
Let $A\to W((t))$ be a $W$-homomorphism, and $A^\dag\to W((t))^\wedge$
the induced homomorphism where $W((t))^\wedge$ denotes the $p$-adic completion.
Suppose that there is an isomorphism $\cX\cong X\times_AW((t))$ and that
$\sigma$ and $\wh\sigma$ are compatible under the map
$A^\dag\to W((t))^\wedge$.
Then the Frobenius $\Phi$ agrees with $\Phi_\crys$ under the natural map
\begin{equation}\label{crys-dR-eq2}
\xymatrix{
H^*_\logcrys((\cY_{\ol\F_p},D_{\ol\F_p})/(W[[t]],(t)))
\ar[r]&  H^*_\dR(X_K/A_K)\ot_AW((t))^\wedge\ar[d]^\cong\\
&
H^*_\rig(X_k/A_k)\ot_{A^\dag}W((t))^\wedge.
}
\end{equation}

\subsection{Explicit description of $\Phi$ by overconvergent functions}\label{exp-formula-sect}
Let 
\[
\xymatrix{
X\ar[r]\ar[d]&Y\ar[d]^f\\
U=\mathrm{Sp}A\ar[r]&\P^1
}
\]
be the fibration in \S \ref{HG-const-sect}.
In what follows we work over the Witt ring $W=W(\ol\F_p)$ with $p>\max(N,M)$.
Put $K:=\Frac W$ the fractional field.

Let $c\in1+pW$ be fixed, and  let 
$\sigma:A^\dag\to A^\dag$ be the $p$-th Frobenius given by
$t^\sigma=ct^p$.
Let
\[
H^1_\rig(X_{\ol\F_p}/A_{\ol\F_p})
\]
be the rigid cohomology group, and $\Phi$ the $\sigma$-linear $p$-th Frobenius.
We shall give an explicit description of $\Phi$.

\begin{lem}\label{free-lem1}
Let $\Spec W[[t]]\to\P^1$, and put $\cY:=Y\times_{\P^1}\Spec W[[t]]$ and
$D:=f^{-1}(0)\subset \cY$ the central fiber. Put $\cX:=\cY\setminus D$.
Then the natural map
\[
H^1(\cY,\Omega^\bullet_{\cY/W[[t]]}(\log D))\lra 
H^1(\cY,\Omega^\bullet_{\cY/W[[t]]}(\log D))\ot W((t))\cong
H^1_\dR(\cX/W((t)))
\]
is injective.
\end{lem}
\begin{pf}
It is enough to show that $H^1(\cY,\Omega^\bullet_{\cY/W[[t]]}(\log D))$ is
$t$-torsion free.
There is an exact sequence
\[
0\lra \vg(\Omega^1_{\cY/W[[t]]}(\log D))
\lra H^1(\cY,\Omega^\bullet_{\cY/W[[t]]}(\log D))\lra H^1(\O_\cY).
\]
The 1st term is $t$-torsion free by Proposition \ref{YdR-prop1}. 
We show that 
$H^1(\O_\cY)$ is a free $W[[t]]$-module.
By \cite[III,12.9]{Ha}, it is enough to show that $\dim_{\kappa(s)}
H^1(Y_s,\O_{Y_s})=(N-1)(M-1)$ for any point $s\in \Spec W[[t]]$
where $\kappa(s)$ is the residue field, and $Y_s:=\cY\times_{W[[t]]}\kappa(s)$.
If $t$ is invertible in $\kappa(s)$, then $Y_s$ is a smooth fiber, and then one has
$\dim_{\kappa(s)}
H^1(Y_s,\O_{Y_s})=(N-1)(M-1)$ as $g(Y_s)=(N-1)(M-1)$.
If $t=0$ in $\kappa(s)$, then $Y_s=D_s:=D\times_W\kappa(s)$ is a simple NCD, 
and then one can directly show that
$\dim_{\kappa(s)}
H^1(D_s,\O_{D_s})=(N-1)(M-1)$.
\end{pf}

The Frobenius $\sigma$ extends on the Frobenius on $K((t))$ as $\sigma(t)=ct^p$.
Let $\Phi_\crys$ be the crystalline Frobenius on 
\[
H^1_\logcrys((\cY_{\ol\F_p},D_{\ol\F_p})/(W[[t]],(t)))
\cong
H^1(\cY,\Omega^\bullet_{\cY/W[[t]]}(\log D)).\] 
Let $1\leq i\leq N-1$ and $1\leq j\leq M-1$ be integers.
Put $a_i:=1-i/N$ and $b_j:=1-j/M$.
Let
\[
F_{ij}(t)=F_{a_i,b_j}(t)={}_2F_1\left({a_i,b_j\atop1};t\right)
\]
be the hypergeometric power series.
It follows from Theorem \ref{YdR-thm1} that the elements
\begin{equation}\label{fermat-form-wt}
\wt\omega_{i,j}:=\frac{1}{F_{ij}(t)}\omega_{i,j},\quad
\wt\eta_{i,j}:=-t(1-t)^{a_i+b_j}F'_{ij}(t)\omega_{ij}+(1-t)^{a_i+b_j-1}F_{ij}(t)\eta_{ij}
\end{equation}
forms a $W[[t]]$-basis of 
\[
H^1(\cY,\Omega^\bullet_{\cY/W[[t]]}(\log D))\cong
\Image[H^1(\Omega^\bullet_{\cY/W[[t]]}(\log D))\to H^1_\dR(\cX/W((t)))]
\]
where the isomorphism follows from Lemma \ref{free-lem1}.

\begin{thm}\label{exp-formula-2}
Let
$\tau_{ij}(t)\in\Q[[t]]$ be defined by
\begin{equation}
\frac{d}{dt}\tau_{ij}(t)=
\frac{1}{t}\left(1-\frac{1}{(1-t)^{a_i+b_j}F_{ij}(t)^2}\right),\quad 
\tau_{ij}(0)=0.
\end{equation}
Let $\psi_p(z)$ be the $p$-adic digamma function introduced in \cite[\S 2]{New}, 
and 
$\log$ the Iwasawa logarithmic function (cf. Appendix A).
Put
\begin{equation}\label{tau}
\tau^{(\sigma)}_{ij}(t)=-2\gamma_p-\psi_p(a_i)-\psi_p(b_j)+p^{-1}\log(c)+\tau_{ij}(t)-p^{-1}\tau_{i'j'}(t^\sigma)\in K[[t]]
\end{equation}
where $i'\in\{1,\ldots,N-1\}$ and $j'\in\{1,\ldots,M-1\}$ are integers 
such that $i'p\equiv i$ mod $N$ and $j'p\equiv j$ mod $M$.
Then
\begin{align}
\label{exp-formula-2-eq1}
\Phi_\crys(\wt\omega_{i'j'})&=p\wt\omega_{ij}+
p\tau^{(\sigma)}_{ij}(t)\wt\eta_{ij}\\
\label{exp-formula-2-eq2}
\Phi_\crys(\wt\eta_{i'j'})&=\wt\eta_{ij}.
\end{align}
\end{thm}

Since
$\Phi_\crys$ agrees with $\Phi$
under the natural map \eqref{crys-dR-eq2},
Theorem \ref{exp-formula-2} implies the following.
\begin{thm}\label{exp-formula-1}
Write $f'(t)=\frac{d}{dt}f(t)$ for a power series $f(t)$.
We define
\begin{align*}
A_{ij}(t)&:=\frac{F_{i'j'}(t^\sigma)}{F_{ij}(t)}-t(1-t)^{a_i+b_j}F'_{ij}(t)
F_{i'j'}(t^\sigma)\tau_{ij}^{(\sigma)}(t)\\
C_{ij}(t)&:=(1-t)^{a_i+b_j-1}F_{ij}(t)F_{i'j'}(t^\sigma)\tau_{ij}^{(\sigma)}(t)\\
B_{ij}(t)&:=pt^\sigma(1-t^\sigma) \frac{F'_{i'j'}(t^\sigma)}{F_{i'j'}(t^\sigma)}A_{ij}(t)
-t\frac{(1-t)^{a_i+b_j}}{(1-t^\sigma)^{a_{i'}+b_{j'}-1}}
\frac{F'_{ij}(t)}{F_{i'j'}(t^\sigma)}\\
D_{ij}(t)&:=
pt^\sigma(1-t^\sigma) \frac{F'_{i'j'}(t^\sigma)}{F_{i'j'}(t^\sigma)}C_{ij}(t)
+\frac{(1-t)^{a_i+b_j-1}}{(1-t^\sigma)^{a_{i'}+b_{j'}-1}}
\frac{F_{ij}(t)}{F_{i'j'}(t^\sigma)}.
\end{align*}
Under the comparison isomorphism
\[
H^1_\rig(X_{\ol\F_p}/A_{\ol\F_p})
\cong H^\bullet_\dR(X/A)\ot_{A} A^\dag_K,
\]
the $p$-th Frobenius $\Phi$ is described as follows,
\[
\begin{pmatrix}
\Phi(\omega_{i'j'})&\Phi(\eta_{i'j'})
\end{pmatrix}
=\begin{pmatrix}
\omega_{ij}&\eta_{ij}
\end{pmatrix}
\begin{pmatrix}
pA_{ij}&B_{ij}\\
pC_{ij}&D_{ij}
\end{pmatrix}.
\]
\end{thm}
\begin{cor}\label{exp-formula-3}
All the power series $\tau^{(\sigma)}_{ij}(t)$, $A_{ij}(t)$, $B_{ij}(t)$, $C_{ij}(t)$ and
$D_{ij}(t)$ lie in the ring $W[[t]]$.
In particular, $A_{ij}(t)$, $B_{ij}(t)$, $C_{ij}(t)$ and
$D_{ij}(t)$ lie in the ring $A^\dag_K\cap W[[t]]=A^\dag\cap W[[t]]$.
\end{cor}
\begin{pf}
Noticing that \eqref{fermat-form-wt} forms a $W[[t]]$-basis, the fact that
$\tau^{(\sigma)}_{ij}(t)\in W[[t]]$ is immediate from Theorem \ref{exp-formula-2} \eqref{exp-formula-2-eq1} together with the fact that
\[\Phi_\crys(\vg(\Omega^1_{\cY/W[[t]]}(\log D))\subset p
H^1(\Omega^\bullet_{\cY/W[[t]]}(\log D)).\]
The others follows from this and the definition.
\end{pf}
\begin{rem}
I don't know a direct proof of Corollary \ref{exp-formula-3} (without $p$-adic cohomology).
\end{rem}
\begin{rem}\label{rem-det}
Note that
$a_{i'}=(a_i)'$ and $b_{j'}=(b_j)'$ (Dwork prime).
In particular $n_i:=a_i-pa_{i'}$ and $m_j:=b_j-pb_{j'}$ are integers $\leq0$.
We have
\begin{equation}\label{det-Frob}
\det\begin{pmatrix}
pA_{ij}&B_{ij}\\
pC_{ij}&D_{ij}
\end{pmatrix}
=p\frac{(1-t)^{a_i+b_j-1}}{(1-t^\sigma)^{a_{i'}+b_{j'}-1}}
=p\frac{(1-t)^{n_i+m_j-1}}{1-t^\sigma}\left(\frac{(1-t)^p}{1-t^\sigma}\right)^{a_{i'}+b_{j'}}
\end{equation}
with
\begin{equation}\label{rem-det-eq1}
\left(\frac{(1-t)^p}{1-t^\sigma}\right)^{a_{i'}+b_{j'}}
=\sum_{n=0}^\infty p^n\binom{a_{i'}+b_{j'}}{n}u(t)^n\in (W[t,(1-t)^{-1}]^\dag)^\times
\end{equation}
where we put $(1-t)^p/(1-t^\sigma)=1+pu(t)$.
In particular
\[
\det\begin{pmatrix}
pA_{ij}&B_{ij}\\
pC_{ij}&D_{ij}
\end{pmatrix}\bigg|_{t=\alpha}=
p\times\text{(unit)}
\]
for $\alpha\in W^\times\setminus(1+pW)$.
\end{rem}
\subsection{Proof of Theorem \ref{exp-formula-2} \eqref{exp-formula-2-eq2}}
For integers $k,l$ with $N {\not{\hspace{-0.2mm}|}}\, k$ and 
$M {\not{\hspace{-0.2mm}|}} \,l$ 
which do not necessarily satisfy that $1\leq k\leq N-1$
and $1\leq l\leq M-1$, $\omega_{kl}$ denotes $\omega_{k_0l_0}$
where $k_0\in \{1,\ldots,N-1\}$ and $l_0\in \{1,\ldots,M-1\}$
such that $k\equiv k_0$ mod $N$ and $l\equiv l_0$ mod $M$.
We apply the same convention to symbols $\eta_{kl}$, $\tau_{kl}(t)$, $a_k$, $b_l$ etc.

\medskip

Let
\[
\nabla:H^1(\cY,\Omega^\bullet_{\cY/W[[t]]}(\log D))\lra
\frac{dt}{t}\ot H^1(\cY,\Omega^\bullet_{\cY/W[[t]]}(\log D))
\]
be the Gauss-Manin connection.
By Proposition \ref{fermat-GM} (or \cite[Prop 4.15]{New}), 
\begin{equation}\label{fermat-GM-wt}
\begin{pmatrix}
\nabla(\wt\omega_{i,j})&\nabla(\wt\eta_{i,j})
\end{pmatrix}
=dt\ot\begin{pmatrix}
\wt\omega_{i,j}&\wt\eta_{i,j}
\end{pmatrix}
\begin{pmatrix}
0&0\\
t^{-1}(1-t)^{-a_i-b_j}F_{ij}(t)^{-2}&0
\end{pmatrix}.
\end{equation}
Using this, one can show
\[
\ker(\nabla)=\bigoplus_{i,j}W\wt\eta_{ij}.
\]
Since $\nabla\Phi_\crys=\Phi_\crys\nabla$, one has
\begin{equation}\label{exp-formula-2-pf-eq1}
\Phi_\crys(\wt\eta_{ij})=\sum_{k,l}\alpha_{kl}\wt\eta_{kl}
\end{equation}
with some constants $\alpha_{kl}\in W$.
Let $i:D\to \cY$ be the embedding.
Let $h$ be the composition as follows
\[
\xymatrix{
\ker(\nabla)\ar[d]\ar[rrd]^h\\
H^1(\Omega^\bullet_{\cY/W[[t]]}(\log D))\ar[r]&H^1(\cY,\O_\cY)\ar[r]_{i^*\quad}
&H^1(D,\O_D).
}
\]
Recall from \S \ref{HG-const-sect} that $D=f^{-1}(0)$ is a simple relative NCD, 
and the irreducible components are $\{D_{x=\zeta_1},D_{y=\zeta_2}
\mid\zeta_1\in \mu_N,\zeta_2\in\mu_M\}$ where 
$D_{x=\zeta_1}:=\{x=\zeta_1\}$ and
$D_{y=\zeta_2}:=\{y=\zeta_2\}$
and $\mu_n:=\{\zeta\in W\mid\zeta^n=1\}$.
Put $P(\zeta_1,\zeta_2):=D_{x=\zeta_1}\cap D_{y=\zeta_2}$ a single point,
and $P:=\{P(\zeta_1,\zeta_2)\}_{\zeta_1,\zeta_2}\subset D$.
There is an exact sequence
\[
\bigoplus_{\zeta_1}H^0(\O_{D_{x=\zeta_1}})
\op\bigoplus_{\zeta_2}H^0(\O_{D_{y=\zeta_2}})
\to
\bigoplus_{\zeta_1,\zeta_2} H^0(\O_{P(\zeta_1,\zeta_2)})
\os{\delta}{\to} 
H^1(\O_D)\to0 
\]
arising from an exact sequence
\[
0\lra \O_D\os{j}{\lra} 
\bigoplus_{\zeta_1}\O_{D_{x=\zeta_1}}
\op\bigoplus_{\zeta_2}\O_{D_{y=\zeta_2}}
\os{u}{\lra}
\bigoplus_{\zeta_1,\zeta_2} \O_{P(\zeta_1,\zeta_2)}\lra 0
\]
where $j$ is the pull-back and $u$ is the map which sends
$(f_{\zeta_1})_{\zeta_1}\times (g_{\zeta_2})_{\zeta_2}$ to
$((g_{\zeta_2}-f_{\zeta_1})|_{P(\zeta_1,\zeta_2)})_{\zeta_1,\zeta_2}$.
\begin{lem}\label{exp-formula-2-pf-lem1}
Let
\[
e_{ij}:=(\zeta_1^i\zeta_2^j)_{\zeta_1,\zeta_2}\in
\bigoplus_{\zeta_1,\zeta_2} H^0(\O_{P(\zeta_1,\zeta_2)})
\]
be an element for $i,j\in \Z$. 
Then $\delta(e_{ij})=h(\wt\eta_{ij})$ for $i\in \{1,\ldots,N-1\}$
and $j\in \{1,\ldots,M-1\}$.
In particular $h\ot\Q$ is bijective.
\end{lem}
\begin{pf}
Recall from \eqref{H01-cech-eq3} the cocycle $(f_{ab,cd})\times(\eta_{ij}^{ab})$
where 
\[
f_{00,11}=f_{00,01}:=x^iy^{j-M},\quad f_{10,11}=f_{10,01}:=(1-t)^2z^{2N-i}y^{j-M}
=(1-t)(1-y^M+z^Ny^M)z^{N-i}y^{j-M},
\]
\[
f_{01,11}:=0,\quad
f_{00,10}:=x^iy^{j-M}-(1-t)^2z^{2N-i}y^{j-M}=(1-x^N)(x^Ny^M-2x^N-y^M)z^{2N-i}y^j.
\]
Note that $D\subset U_{00}\cup U_{11}$.
We have $h(\wt\eta_{ij})=[(x^iy^{j-M}|_D)]$
under the isomorphism
\[
H^1(\O_D)\cong \Coker
\left[\bigoplus_{(a,b)=(0,0),(1,1)}\vg(U_{ab},\O_D)\os{d}{\lra}\vg(U_{00,11},\O_D)
\right]
\]
where $U_{00,11}:=U_{00}\cap U_{11}$ and $d(f_{00},f_{11}):=f_{11}-f_{00}$.
A diagram chase 
\[
\xymatrix{
&\bigoplus_{\zeta_1}\vg(U_{ab},\O_{D_{x=\zeta_1}})
\times\bigoplus_{\zeta_2}\vg(U_{ab},\O_{D_{y=\zeta_2}})
\ar[r]^{\hspace{2cm}u}\ar[d]^d
&\bigoplus_{\zeta_1,\zeta_2}\vg(\O_{P(\zeta_1,\zeta_2)})\\
\vg(U_{00,11},\O_D)\ar[r]&
\bigoplus_{\zeta_1}\vg(U_{00,11},\O_{D_{x=\zeta_1}})\times 
\bigoplus_{\zeta_2}\vg(U_{00,11},\O_{D_{y=\zeta_2}})
}
\]
\[
\xymatrix{
&(0,\zeta_1^iw^{M-j})
\times(-\zeta_2^jx^i,0)\ar[r]^{\hspace{1cm}u}
\ar[d]^d&(\zeta_1^i\zeta_2^j)=e_{ij}\\
h(\wt\eta_{ij})=(x^iy^{j-M}|_D)\ar[r]&(\zeta_1^iy^{j-M},\zeta_2^jx^i)
}
\]
yields $\delta(e_{ij})=h(\wt\eta_{ij})$. The last statement is an exercise of linear algebra.
\end{pf}
We turn to the proof of \eqref{exp-formula-2-eq2}.
Apply $\Phi_\crys$ on the equality $h(\wt\eta_{ij})=\delta(e_{ij})$
in Lemma \ref{exp-formula-2-pf-lem1}.
Since $h$ and $\delta$ are compatible with respect to the action of $\Phi_\crys$,
one has
\[
h\Phi(\wt\eta_{ij})=\sum_{k,l}\alpha_{kl}h(\wt\eta_{kl})=\delta\Phi_\crys(e_{ij})
\]
by \eqref{exp-formula-2-pf-eq1}.
On the other hand
\[
\Phi_\crys(e_{ij})=
(\zeta_1^{ip}\zeta_2^{jp})_{\zeta_1,\zeta_2}=e_{ip,jp}\in
\bigoplus_{\zeta_1,\zeta_2} H^0(\O_{P(\zeta_1,\zeta_2)})
\]
by definition of $\Phi_\crys$.
Therefore one has
\[
\sum_{k,l}\alpha_{kl}h(\wt\eta_{kl})=\delta(e_{ip,jp})=h(\wt\eta_{ip,jp}),
\]
and hence $\alpha_{kl}=1$ if $(k,l)=(ip,jp)$ in $\Z/N\Z\times \Z/M\Z$
and $=0$ otherwise.
This completes the proof of \eqref{exp-formula-2-eq2}.
\subsection{Proof of Theorem \ref{exp-formula-2} \eqref{exp-formula-2-eq1}}
For $(\zeta_1,\zeta_2)\in \mu_N\times\mu_M$, we denote by $[\zeta_1,\zeta_2]$ 
the automorphism of $\cY$ given by $(x,y,t)\mapsto(\zeta_1x,\zeta_2y,t)$.
Since $[\zeta_1,\zeta_2]\Phi_\crys=\Phi_\crys[\zeta_1,\zeta_2]$, one has
\[
\Phi_\crys(\wt\omega_{ij})\in W[[t]]\wt\omega_{ip,jp}+W[[t]]\wt\eta_{ip,jp}.
\]
One can further show that there is $g_{ij}(t)\in W[[t]]$ such that
\begin{equation}\label{exp-formula-2-pf-eq2}
\Phi_\crys(\wt\omega_{ij})=p\omega_{ip,jp}+g_{ij}(t)\wt\eta_{ip,jp}
\end{equation}
(this can be proved in the same way as the proof of \cite[Lemma 4.5]{New}).
Thus our goal is to show $g_{ij}(t)=p\tau^{(\sigma)}_{ip,jp}(t)$.
Apply $\nabla$ on \eqref{exp-formula-2-pf-eq2}.
It follows from \eqref{fermat-GM-wt} that we have
\begin{align*}
\text{LHS}&=\nabla\Phi_\crys(\wt\omega_{ij})\\
&=\Phi_\crys\nabla(\wt\omega_{ij})\\
&=
\Phi_\crys\left((1-t)^{-a_i-b_j}F_{ij}(t)^{-2}
\frac{dt}{t}\ot\wt\eta_{ij}\right)\\
&=
p(1-t^\sigma)^{-a_i-b_j}F_{ij}(t^\sigma)^{-2}
\frac{dt}{t}\ot\wt\eta_{ip,jp}\quad \text{(by Theorem \ref{exp-formula-2}
\eqref{exp-formula-2-eq2})}
\end{align*}
and
\[
\text{RHS}=p(1-t)^{-a_{ip}-b_{jp}}F_{ip,jp}(t)^{-2}
\frac{dt}{t}\ot\wt\eta_{ip,jp}+g'_{ij}(t)dt\ot\wt\eta_{ip,jp}.
\]
Hence
\[
g'_{ij}(t)=\frac{p}{t}\left[\frac{1}{(1-t^\sigma)^{a_i+b_j}F_{ij}(t^\sigma)^2}
-\frac{1}{(1-t)^{a_{ip}+b_{jp}}F_{ip,jp}(t)^2}\right]
\]
or equivalently
\begin{equation}\label{exp-formula-2-pf-eq3}
g_{ij}(t)=p(C_{ij}+\tau_{ip,jp}(t)-p^{-1}\tau_{ij}(t^\sigma))
\end{equation}
with $C_{ij}$ a constant.
The rest is to show 
\begin{equation}\label{exp-formula-2-pf-eq4}
C_{ij}=-2\gamma_p-\psi_p(a_{ip})-\psi_p(b_{jp})+p^{-1}\log(c).
\end{equation}
To do this, we recall from \cite[4.6]{New}
the {\it regulator formula}.

\medskip

For $(\nu_1,\nu_2)\in\mu_N(K)\times \mu_M(K)$,
let
\begin{equation}\label{m-fermat-eq1}
\xi=\xi(\nu_1,\nu_2)=\left\{
\frac{x-1}{x-\nu_1},\frac{y-1}{y-\nu_2}
\right\}\in K_2(X)
\end{equation}
be a $K_2$-symbol.
The symbol $\xi$ defines the $1$-extension
\[
0\lra H^1(X/A)(2)\lra M_\xi(X/A)\lra A\lra 0
\]
in the category of $\FilFMIC(A)$ (see \cite[4.5]{New} or \cite[2.5]{AM} 
for the notation).
Let $e_\xi\in \mathrm{Fil}^0M_\xi(X/A)_\dR$ be the unique lifting of $1\in A$.
Let $E_k^{(ij)}(t)\in W[[t]]$ be defined by
\begin{equation}\label{fermat-e-eq2}
e_\xi-\Phi_\crys(e_\xi)
=-N^{-1}M^{-1}
\sum_{i=1}^{N-1}\sum_{j=1}^{M-1}(1-\nu^{-i}_1)(1-\nu^{-j}_2)[E^{(ij)}_1(t)\wt\omega_{ij}
+E^{(ij)}_2(t)\wt\eta_{ij}].
\end{equation}
Then one of the main results in \cite{New} is 
\begin{equation}\label{exp-formula-2-pf-eq5}
\frac{E_1^{(ij)}(t)}{F_{ij}(t)}=-\cF^{(\sigma)}_{a_i,b_j}(t)
\end{equation}
(\cite[Theorem 4.18]{New}) where $\cF^{(\sigma)}_{\ul a}(t)$ is the 
$p$-adic hypergeometric function of log type introduced in \cite[\S 3]{New}.

\medskip

We turn to the proof of \eqref{exp-formula-2-pf-eq4}.
Apply $\nabla$ on \eqref{fermat-e-eq2}.
Noticing that $\Phi_\crys\nabla=\nabla\Phi_\crys$ and
\[
\nabla(e_\xi)=\mathrm{dlog}(\xi)=N^{-1}M^{-1}
\sum_{i=1}^{N-1}\sum_{j=1}^{M-1}(1-\nu^{-i}_1)(1-\nu^{-j}_2)\frac{dt}{t}
\omega_{ij},
\]
one has a differential equation
\[
t\frac{d}{dt}E^{(ij)}_2(t)+(1-t)^{-a_i-b_j}F_{ij}(t)^{-2}E^{(ij)}_1(t)
=p^{-1}F_{i'j'}(t^\sigma)g_{i'j'}(t^\sigma)
\]
by \eqref{fermat-GM-wt} and \eqref{exp-formula-2-pf-eq2}
where $i',j'$ are integers such that $i'\in \{1,\ldots,N-1\}$ with $pi'\equiv i$ mod $N$
and 
$j'\in \{1,\ldots,M-1\}$ with $pj'\equiv j$ mod $M$.
Substitute $t=0$ in the above. We have
\[E^{(ij)}_1(0)=p^{-1}g_{i'j'}(0)=C_{i'j'}.\]
By \eqref{exp-formula-2-pf-eq5}, 
\[
E^{(ij)}_1(0)=-\cF^{(\sigma)}_{a_{i'},b_{j'}}(0)=
-2\gamma_p-\psi_p(a_{i'})-\psi_p(b_{j'})+p^{-1}\log(c),
\]
and hence \eqref{exp-formula-2-pf-eq4}
as required. This completes the proof of Theorem \ref{exp-formula-2} \eqref{exp-formula-2-eq1}.

\section{Computing Dwork's $p$-adic Hypergeometric functions}
In this section, we shall give an algorithm for computing special values
of Dwork's $p$-adic hypergeometric functions whose bit complexity increases at most 
$O(n^4\log^3n)$ as $n\to\infty$.
\subsection{$p$-adic expansions of $A_{ij}(t),B_{ij}(t),C_{ij}(t),D_{ij}(t)$}
We keep the setting in \S \ref{exp-formula-sect}.
Recall Theorem \ref{exp-formula-1},
\[
\begin{pmatrix}
\Phi(\omega_{i'j'})&\Phi(\eta_{i'j'})
\end{pmatrix}
=\begin{pmatrix}
\omega_{ij}&\eta_{ij}
\end{pmatrix}
\begin{pmatrix}
pA_{ij}(t)&B_{ij}(t)\\
pC_{ij}(t)&D_{ij}(t)
\end{pmatrix}
\]
with 
$A_{ij}(t),B_{ij}(t),C_{ij}(t),D_{ij}(t)\in A^\dag\cap W[[t]]$ (Corollary \ref{exp-formula-1}).
By \cite[Theorem 2.1]{KT}, the overconvergent functions
$A_{ij}(t),\ldots,D_{ij}(t)$ have `nice' $p$-adic expansions, and this is the key fact in
our algorithm.
We here write down the necessary statement.
 \begin{thm}\label{main-thm}
For an integer $n\geq 1$, define
\[e_n:=\max\{k\in \Z_{\geq 1}\mid \ord_p(p^k/k!)<n\}.\]
Then
\[
pA_{ij}(t)\equiv p\frac{(\text{polynomial of degree}\leq pe_n+p)}{(1-t^\sigma)(1-t)^{pe_n}}\mod p^nW[[t]],
\]
\[
B_{ij}(t)\equiv \frac{(\text{polynomial of degree}\leq pe_n+2p)}{(1-t^\sigma)(1-t)^{pe_n}}\mod p^nW[[t]],
\]
\[
pC_{ij}(t)\equiv p\frac{(\text{polynomial of degree}\leq pe_n+p-1)}{(1-t^\sigma)(1-t)^{pe_n}}\mod p^nW[[t]],
\]
\[
D_{ij}(t)\equiv \frac{(\text{polynomial of degree}\leq pe_n+2p-1)}
{(1-t^\sigma)(1-t)^{pe_n}}\mod p^nW[[t]].
\]
\end{thm}
\begin{rem}
Since $p\ne 2$ by the assumption, 
$e_n<\infty$ for any $n\geq1$. More precisely
\[
e_n\sim \frac{p-1}{p-2}n\quad \text{as }n\to \infty.
\]
\end{rem}
\begin{rem}
The degrees $pe_n+p$ etc. are not optimal.
\end{rem}
We give a self-contained proof of Theorem \ref{main-thm} for the sake of the completeness.
\begin{pf} (cf. \cite[p.11--13]{KT}).
Let $\l=t$, $1-t$ or $t^{-1}$.
Let $\sigma_\l$ be the $p$-th Frobenius on $W[[\l]]$ given by $\sigma_\l(\l)=\l^p$.
Note that $\sigma_\l$ induces the $p$-th Frobenius on 
$A^\dag=W[t,(t-t^2)^{-1}]^\dag$.
Let $\Phi_\l$ denote the $\sigma_\l$-linear Frobenius on 
\[
H^1_\rig(X_{\ol\F_p}/A_{\ol\F_p})
\cong H^1_\dR(X/A)\ot_{A} A^\dag_K.
\]
Let $\sigma$ be the Frobenius given by $\sigma(t)=ct^p$ and 
$\Phi$ the $\sigma$-linear Frobenius as in \S \ref{exp-formula-sect}.
Then the relation with $\Phi_\l$ is given as follows (\cite[6.1]{EK}, \cite[17.3.1]{Ke}).
\begin{equation}\label{f-trans}
\Phi(x)-\Phi_\l(x)=\sum_{k=1}^\infty\frac{(\l^\sigma-\l^p)^k}{k!}\Phi_\l\partial^k_\l x,
\quad  x\in H^1_\dR(X/A)\ot_{A} A^\dag_K
\end{equation}
where $\partial_\l:=\nabla_{d/d\l}$.
Let $\l=1-t$.
Since $\l^\sigma-\l^p=pw(\l)\in pW[\l]$,
\eqref{f-trans} yields
\[
\Phi(x)-\Phi_\l(x)=\sum_{k=1}^\infty\frac{p^k}{k!}w(t)^k\Phi_\l\partial^k_\l x.
\]
Note that $\Phi_\l(H_\l)\subset H_\l$ while $\Phi(H_\l)\not\subset H_\l$.
Since 
$\partial^k_\l(H_\l)\subset \l^{-k}H_\l$, one has
$\Phi_\l\partial^k_\l(H_\l)\subset \l^{-kp}H_\l$ for all $k\geq0$, and hence
\[
\Phi(x)\in \sum_{k=0}^\infty\frac{p^k}{k!}\l^{-kp}H_\l.
\]
We thus have
\begin{equation}\label{mt-pf-eq1}
\Phi(x)\in\l^{-pe_n}H_\l+p^n\wh H_\l,\quad\forall\, n\geq 1,
\,\forall x\in (H^1_\dR(X/A)\ot_{A} A^\dag_K)\cap H_\l
\end{equation}
if $\l=1-t$
where $\wh H_\l$ is the $p$-adic completion of $H_\l\ot_{W[[\l]]} W((\l))$.
Let $\l=t^{-1}$.
In this case, since $\sigma_\l(t)=t^p$, the Frobenius $\Phi_\l$ acts on the $W[[\l]]$-lattice $H_\l$. Hence 
\begin{equation}\label{mt-pf-eq2}
\Phi_\l(H_\l)\subset H_\l,\quad \l=t^{-1}.
\end{equation}

Let us prove Theorem \ref{main-thm}.
Since $A_{ij},B_{ij},C_{ij},D_{ij}\in A^\dag\cap W[[t]]$, one can write
\[
pA_{ij}(t)\text{ mod }p^nW[[t]]= \frac{pF^A_{ij}(t)}{(1-t^\sigma)(1-t)^{d^A_{ij}}},
\]
\[
B_{ij}(t)\text{ mod }p^nW[[t]]= \frac{F^B_{ij}(t)}{(1-t^\sigma)(1-t)^{d^B_{ij}}},
\]
\[
pC_{ij}(t)\text{ mod }p^nW[[t]]= \frac{pF^C_{ij}(t)}{(1-t^\sigma)(1-t)^{d^C_{ij}}},
\]
\[
D_{ij}(t)\text{ mod }p^nW[[t]]= \frac{F^D_{ij}(t)}{(1-t^\sigma)(1-t)^{d^D_{ij}}},
\]
in $W/p^nW[[t]]$ with $F^A_{ij}(t),\, F^B_{ij}(t),\ldots\in W/p^nW[t]$ polynomials and 
$d_{ij}^A,\, d^B_{ij},\ldots\in \Z_{\geq0}$.
Let $\l=t^{-1}$. Then $H_\l$ is a free $W[[t]]$-module with
basis $\{\omega_{ij},\l\eta_{ij}\}$ (Theorem \ref{YdR-thm1} (2)).
Therefore 
it follows from \eqref{mt-pf-eq2} that the entries of the $2\times 2$-matrix in below
lie in $W[[\l]]$,
\[
\begin{pmatrix}
\Phi(\omega_{i'j'})&\Phi(\l\eta_{i'j'})
\end{pmatrix}
=\begin{pmatrix}
\omega_{ij}&\l\eta_{ij}
\end{pmatrix}
\begin{pmatrix}
pA_{ij}&\l^\sigma B_{ij}\\
p\l^{-1}C_{ij}&\l^\sigma\l^{-1}D_{ij}
\end{pmatrix}.
\]
This implies
\begin{equation}\label{mt-pf-eq4}
\begin{cases}
\deg(pF^A_{ij})\leq d_{ij}^A+p\\
\deg(F^B_{ij})\leq d_{ij}^B+2p\\
\deg(pF^C_{ij})\leq d_{ij}^C+p-1\\
\deg(F^D_{ij})\leq d_{ij}^D+2p-1.
\end{cases}
\end{equation}
Next we give upper bounds of $d^A_{ij},d^B_{ij},d^C_{ij}$ and $d^D_{ij}$.
Let $\l=1-t$ and 
let $\omega^*_{ij},\eta_{ij}^*$ be the basis of
$H_\l$ in Theorem \ref{YdR-thm1} (3).
Let
\[
\begin{pmatrix}
\Phi(\omega^*_{i'j'})&\Phi(\eta^*_{i'j'})
\end{pmatrix}
=\begin{pmatrix}
\omega^*_{ij}&\eta^*_{ij}
\end{pmatrix}
\begin{pmatrix}
pA^*_{ij}&B^*_{ij}\\
pC^*_{ij}&D^*_{ij}
\end{pmatrix}.
\]
It follows from \eqref{mt-pf-eq1} that we have
\begin{equation}\label{mt-pf-eq5}
(1-t)^{pe_n}pA^*_{ij},\,
(1-t)^{pe_n}B^*_{ij},\,
(1-t)^{pe_n}pC^*_{ij},\,
(1-t)^{pe_n}D^*_{ij}\in W[[\l]]+p^nW((\l))^\wedge.
\end{equation}
If $i/N+j/M\geq 1$, then $(\omega^*_{ij},\eta^*_{ij})=(\omega_{ij},\eta_{ij})$,
and if $i/N+j/M> 1$, then
\[
\begin{pmatrix}
\omega^*_{ij}&\eta^*_{ij}
\end{pmatrix}
=\begin{pmatrix}
\omega_{ij}&\eta_{ij}
\end{pmatrix}
\begin{pmatrix}
1-t&lt\\
0&-1
\end{pmatrix}
\]
where $l:=1-i/N-j/M$.
Therefore if $i/N+j/M\geq 1$ and $i'/N+j'/M\geq 1$, then
\[
\begin{pmatrix}
pA_{ij}&B_{ij}\\
pC_{ij}&D_{ij}
\end{pmatrix}
=\begin{pmatrix}
pA^*_{ij}&B^*_{ij}\\
pC^*_{ij}&D^*_{ij}
\end{pmatrix}.
\]
By \eqref{mt-pf-eq5}, we have 
$d_{ij}^A,\,
d_{ij}^B,\,
d_{ij}^C,\,
d_{ij}^D\leq pe_n$.
If $i/N+j/M< 1$ and $i'/N+j'/M\geq 1$, then
\[
\begin{pmatrix}
pA_{ij}&B_{ij}\\
pC_{ij}&D_{ij}
\end{pmatrix}
=\begin{pmatrix}
1-t&lt\\
0&-1
\end{pmatrix}\begin{pmatrix}
pA^*_{ij}&B^*_{ij}\\
pC^*_{ij}&D^*_{ij}
\end{pmatrix}.
\]
By \eqref{mt-pf-eq5}, we have 
$d_{ij}^A\leq pe_n-1$ and $
d_{ij}^B,\,
d_{ij}^C,\,
d_{ij}^D\leq pe_n$.
If $i/N+j/M\geq 1$ and $i'/N+j'/M< 1$, then
\begin{align*}
\begin{pmatrix}
pA_{ij}&B_{ij}\\
pC_{ij}&D_{ij}
\end{pmatrix}
&=\begin{pmatrix}
pA^*_{ij}&B^*_{ij}\\
pC^*_{ij}&D^*_{ij}
\end{pmatrix}
\begin{pmatrix}
1-t^\sigma&lt^\sigma\\
0&-1
\end{pmatrix}^{-1}\\
&=\frac{1}{t^\sigma-1}\begin{pmatrix}
-pA^*_{ij}&(1-t^\sigma)B^*_{ij}-lt^\sigma pA_{ij}^*\\
-pC^*_{ij}&(1-t^\sigma)D^*_{ij}-lt^\sigma pC_{ij}^*
\end{pmatrix}.
\end{align*}
By \eqref{mt-pf-eq5}, we have 
$d_{ij}^A,\,
d_{ij}^B,\,
d_{ij}^C,\,
d_{ij}^D\leq pe_n$.
If $i/N+j/M<1$ and $i'/N+j'/M< 1$, then
\begin{align*}
\begin{pmatrix}
pA_{ij}&B_{ij}\\
pC_{ij}&D_{ij}
\end{pmatrix}
&=\begin{pmatrix}
1-t&lt\\
0&-1
\end{pmatrix}
\begin{pmatrix}
pA^*_{ij}&B^*_{ij}\\
pC^*_{ij}&D^*_{ij}
\end{pmatrix}
\begin{pmatrix}
1-t^\sigma&lt^\sigma\\
0&-1
\end{pmatrix}^{-1}\\
&=\frac{1}{t^\sigma-1}
\begin{pmatrix}
(t-1)pA_{ij}^*-ltpC^*_{ij}&\cdots\\
pC^*_{ij}&(t^\sigma-1)D_{ij}^*+lt^\sigma pC^*_{ij}
\end{pmatrix}.
\end{align*}
By \eqref{mt-pf-eq5}, we have 
$d_{ij}^A,\,
d_{ij}^B,\,
d_{ij}^C,\,
d_{ij}^D\leq pe_n$.
In any case one has
\begin{equation}\label{mt-pf-eq6}
d_{ij}^A,\,
d_{ij}^B,\,
d_{ij}^C,\,
d_{ij}^D\leq pe_n.
\end{equation}
Theorem \ref{main-thm} follows from \eqref{mt-pf-eq4} and \eqref{mt-pf-eq6}.
\end{pf}
\subsection{Algorithm for computing Dwork's $p$-adic hypergeometric functions}
\label{algorithm-sect}
For $a,b\in \Q$, let
\[
F_{ab}(t)={}_2F_1\left({a,b\atop 1};t\right)
\]
be the hypergeometric power series.
We give an algorithm for computing the special values
\begin{equation}\label{algorithm-eq0}
\cF^{\Dw,\sigma}_{ab}(t):=\frac{F_{ab}(t)}{F_{a'b'}(t^\sigma)}
,\quad\frac{F'_{ab}(t)}{F_{ab}(t)}
\end{equation}
at $\alpha\in W^\times\setminus(1+pW)$ modulo $p^n$.

\medskip

\noindent{\bf Notation.}
Let $N,M\geq 2$ be integers, and 
$p>\max(N,M)$ a prime.
Let $W=W(\ol\F_p)$.
Let $a,b\in \Q$ satisfy that $a\in \frac{1}{N}\Z$ and $b\in \frac{1}{M}\Z$
and $0<a,b<1$.
Let $a'$ denote the Dwork prime (see \S \ref{DefHG-sect}).
Let $c\in 1+pW$, and let $\sigma:W[[t]]\to W[[t]]$ be the $p$-th Frobenius
given by $\sigma(t)=ct^p$.
Following the notation in \eqref{tau} and Theorem \ref{exp-formula-1},
we define
\[
\frac{d}{dt}\tau_{ab}(t)=
\frac{1}{t}\left(1-\frac{1}{(1-t)^{a+b}F_{ab}(t)^2}\right),\quad 
\tau_{ab}(0)=0,
\]
\[
\tau^{(\sigma)}_{ab}(t)=-2\gamma_p-\psi_p(a)-\psi_p(b)+p^{-1}\log(c)+\tau_{ab}(t)-p^{-1}\tau_{a'b'}(t^\sigma)\in W[[t]],
\]
and 
\begin{align*}
A_{ab}(t)&:=\frac{F_{a'b'}(t^\sigma)}{F_{ab}(t)}-t(1-t)^{a+b}F'_{ab}(t)
F_{a'b'}(t^\sigma)\tau_{ab}^{(\sigma)}(t)\\
C_{ab}(t)&:=(1-t)^{a+b-1}F_{ab}(t)F_{a'b'}(t^\sigma)\tau_{ab}^{(\sigma)}(t)\\
B_{ab}(t)&:=pt^\sigma(1-t^\sigma) \frac{F'_{a'b'}(t^\sigma)}{F_{a'b'}(t^\sigma)}A_{ab}(t)
-t\frac{(1-t)^{a+b}}{(1-t^\sigma)^{a'+b'-1}}
\frac{F'_{ab}(t)}{F_{a'b'}(t^\sigma)}\\
D_{ab}(t)&:=
pt^\sigma(1-t^\sigma) \frac{F'_{a'b'}(t^\sigma)}{F_{a'b'}(t^\sigma)}C_{ab}(t)
+\frac{(1-t)^{a+b-1}}{(1-t^\sigma)^{a'+b'-1}}
\frac{F_{ab}(t)}{F_{a'b'}(t^\sigma)}.
\end{align*}
Let $a^{(k)}$ denote the $k$-th Dwork prime.
Put
\[
F^{(k)}(t):=F_{a^{(k)}b^{(k)}}(t),\quad 
A^{(k)}_\sigma(t):=A_{a^{(k)}b^{(k)}}(t),\ldots,
D^{(k)}_\sigma(t):=D_{a^{(k)}b^{(k)}}(t),
\]
\[
\DF^{(k)}(t):=\frac{(F^{(k)}(t))'}{F^{(k)}(t)}.
\]
for $k\geq 0$. 
Note that $F^{(k)}(t)$ and $\DF^{(k)}(t)$
do not depend on $\sigma$.
We put
\[
E^{(k)}_\sigma(t):
=\frac{(1-t)^{a^{(k)}+b^{(k)}}}{(1-t^\sigma)^{a^{(k+1)}+b^{(k+1)}}}
=(1-t)^{m_k}\left(\frac{(1-t)^p}{1-t^\sigma}\right)^{a^{(k+1)}+b^{(k+1)}}
\]
where $m_k:=a^{(k)}-pa^{(k+1)}+b^{(k)}-pb^{(k+1)}\in \Z_{\leq 0}$.
Note $E^{(k)}_\sigma(t)\in (A^\dag)^\times$.
We have
\begin{equation}\label{algorithm-eq1}
D^{(k)}_\sigma(t)=
pt^\sigma(1-t^\sigma)C_\sigma^{(k)}(t)\DF^{(k+1)}(t^\sigma)+
\frac{1-t^\sigma}{1-t}E^{(k)}_\sigma(t)
\cF^{\Dw,(k)}_\sigma(t)
\end{equation}
and
\begin{align}
&\overbrace{\begin{pmatrix}
pA_\sigma^{(k)}(t)&B_\sigma^{(k)}(t)\\
pC_\sigma^{(k)}(t)&D_\sigma^{(k)}(t)
\end{pmatrix}}^{H_\sigma^{(k)}(t)}
\begin{pmatrix}
t^\sigma(1-t^\sigma)\DF^{(k+1)}(t^\sigma)\\
-1
\end{pmatrix}\notag\\
&\qquad=\frac{1-t^\sigma}{1-t}E^{(k)}_\sigma(t)\cF^{\Dw,\sigma}_{a^{(k)}b^{(k)}}(t)
\begin{pmatrix}
t(1-t)\DF^{(k)}(t)\\
-1
\end{pmatrix}.\label{algorithm-eq2}
\end{align}

\medskip
\begin{center}
{\large\bf Algorithm}
\end{center}

Let $m\geq1$ be the smallest integer such that $(a^{(m)},b^{(m)})=(a,b)$.
Let $\alpha\in W^\times\setminus(1+pW)$ be an arbitrary element satisfying
\[
F^{(k)}(t)_{<p}|_{t=\alpha}\not\equiv 0\mod pW,\quad 0\leq \forall\,k\leq m-1.
\] 
Let $\sigma(t)=ct^p$ with $c\in 1+pW$ arbitrary.
The algorithm for computing \eqref{algorithm-eq0} is the following.

\medskip

\noindent{\bf Step 1}.
Let $\beta\in W^\times\setminus(1+pW)$
satisfy
\[
F^{(k)}(t)_{<p}|_{t=\beta}\not\equiv 0\mod pW,\quad 0\leq \forall\,k\leq m-1.
\] 
In {\bf Step 3}, we shall take $\beta=t^\sigma|_{t=\alpha}=c\alpha^p$.
Let $\sigma_\beta(t)=\beta^{1-p}t^p$ so that we have $t^{\sigma_\beta}
|_{t=\beta}=\beta$.
Then we compute the special values
\[
pA^{(k)}_{\sigma_\beta}(\beta),\,
pC^{(k)}_{\sigma_\beta}(\beta),\,
B^{(k)}_{\sigma_\beta}(\beta),\,
D^{(k)}_{\sigma_\beta}(\beta)
\mod p^nW
\]
for each $k=0,1,\ldots,m-1$.
One can do it in the following way.
Compute the power series
\[
(1-t^{\sigma_\beta})(1-t)^{pe_n}A^{(k)}_{\sigma_\beta}(t)
\]
until the degree $pe_n+p$, say $F^A(t)$. Then
it follows from Theorem \ref{main-thm} that
\[
pA^{(k)}_{\sigma_\beta}(t)\equiv \frac{pF^A(t)}{(1-t^{\sigma_\beta})(1-t)^{pe_n}}
\mod p^nW[[t]
\]
and hence
\[
pA^{(k)}_{\sigma_\beta}(\beta)\equiv
\frac{pF^A(\beta)}{(1-\beta)^{pe_n+1}}\mod p^nW.
\]
The other values are obtained in the same way.

\medskip

\noindent{\bf Step 2}.
We mean $H^{(l)}_{\sigma_\beta}=H^{(l_0)}_{\sigma_\beta}$ for arbitrary $l\in \Z$
where $l_0\in\{0,1,\ldots,m-1\}$
such that $l\equiv l_0$ mod $m$.

Compute an eigenvector ${\mathbf u}_\beta$ of a $2\times 2$-matrix
\[
H^{(k-m)}_{\sigma_\beta}(\beta)\cdots H^{(k-2)}_{\sigma_\beta}(\beta)
H^{(k-1)}_{\sigma_\beta}(\beta)
\]
whose eigenvalue is a unit.
This is unique up to scalar.
Indeed, it follows from \eqref{algorithm-eq2}
that the vector
\begin{equation}\label{algorithm-eq3}
\begin{pmatrix}
\beta(1-\beta)\DF^{(k)}(\beta)\\
-1
\end{pmatrix}
\end{equation}
is an eigenvector of 
$H^{(1)}_{\sigma_\beta}(\beta)\cdots H^{(m)}_{\sigma_\beta}(\beta)$
whose eigenvalue is
\begin{equation}\label{algorithm-eq4}
\prod_{k=0}^{m-1}E^{(k)}_{\sigma_\beta}(\beta)
\cF^{\Dw,\sigma_\beta}_{a^{(k)}b^{(k)}}(\beta)\in W^\times.
\end{equation}
The other eigenvalue is not a unit as
$\det(H^{(0)}_{\sigma_\beta}(\beta)\cdots 
H^{(m-1)}_{\sigma_\beta}(\beta))=p^m\times$(unit) by Remark \ref{det-Frob}
(actually the determinant is equal to $p^m$).
Therefore \eqref{algorithm-eq3} 
is characterized as the eigenvector with the unique eigenvalue which is a unit.
We thus have the special value
\[
\DF^{(k)}(\beta)=\frac{F'_{a^{(k)}b^{(k)}}(t)}{F_{a^{(k)}b^{(k)}}(t)}\bigg|_{t=\beta}
\mod p^nW
\]
for each $k$.

\medskip

\noindent{\bf Step 3}.
Let $\sigma(t)=ct^p$ be as in the beginning.
Take $\beta=t^\sigma|_{t=\alpha}=c\alpha^p$ in {\bf Step 2}.
We have
\[
\DF^{(1)}(\beta)=
\DF^{(1)}(t^\sigma)|_{t=\alpha}
\mod p^nW.
\]
Compute the special values
\[
pC^{(0)}_\sigma(\alpha),\,
D^{(0)}_\sigma(\alpha),\,
E^{(0)}_\sigma(\alpha)\mod p^nW
\]
according to {\bf Step 1}, and
\[
E^{(0)}_\sigma(\alpha)\mod p^nW
\]
utilizing the expansion
\[
\left(\frac{(1-t)^p}{1-t^\sigma}\right)^{a+b}
=\sum_{n=0}^\infty p^n\binom{a+b}{n}u(t)^n,\quad
\frac{(1-t)^p}{1-t^\sigma}=1+pu(t).
\]
Substitute $t=\alpha$ in \eqref{algorithm-eq1}.
Then we have the special value
\[
\cF^{\Dw,\sigma}_{ab}(\alpha)
=\frac{F_{ab}(t)}{F_{a'b'}(t^\sigma)}\bigg|_{t=\alpha}\mod p^nW
\]
as $E^{(0)}_\sigma(\alpha)\in W^\times$.

\subsection{Bit Complexity}
We give an upper estimate of the bit complexity of the algorithm 
displayed in \S \ref{algorithm-sect}.

\medskip

We review the notion of the bit complexity. A general reference is
 the text book \cite{Zimmermann}.
The bit of a natural number $N$ is defined to be the number of digits 
of $N$ in binary notation, so it is at most $\log_2 (N+1)$.
The bit of $N!$ is at most $\log_2(N!+1)\sim (\log 2)^{-1}N\log N$ (Stirling).
The bit complexity of an algorithm is defined to be the number of single operations
to complete the algorithm.  The bit complexity of
(1-digit)$\pm$(1-digit) or (1-digit)$\times$(1-digit) is 1 by definition.
We denote by $M(n,m)$ the bit complexity of multiplication
($n$-digits$)\times(m$-digits).
We write $M(n)=M(n,n)$.
By the naive multiplication algorithm,
$M(n,m)$ is $O(nm)$,
which means that there is a constant $C$ such that 
$M(n,m)\leq Cnm$ when $n,m\to \infty$.
We sum up the basic results.
\begin{itemize}
\item
For integers $i,j\geq 0$, 
the bit complexity of $i\pm j$ is $O(\max(\log i,\log j))$.
\item
The bit complexity of $i\cdot j$ is $M(\log i,\log j)$ (which is at most 
$O(\log i\log j)$).
\item
The bit complexity for computing the remainder ($i$ mod $j$) is 
$M(\log i,\log j)$.
\end{itemize}
Let $a$ be a fixed rational number. Then
the bit complexity of $(a)_i$ is at most 
\begin{equation}\label{bit-eq1}
\sum_{n=1}^iM(n\log n,\log n)\leq O(i^2(\log i)^2)
\end{equation}
by computing it in the following way
\[
(a)_i=(a+i-1)\cdot(a)_{i-1},\quad
(a)_{i-1}=(a+i-2)\cdot(a)_{i-2},\ldots
\]
Let $a_i,b_j$ be rational numbers whose denominators and numerators
are less than $k$.
Let $f(t)=\sum_{i=0}^na_it^i$ and $g(t)=\sum_{j=0}^nb_jt^j$.
Then the bit complexity of computing $f(t)\pm g(t)$ 
is $O(n\log k)$
The bit complexity of computing $f(t)g(t)$
\begin{equation}\label{bit-eq2}
n^2M(\log k)+O(n^2\log (kn))<O(n^2(\log n+(\log k)^2))
\end{equation}
on noticing that the coefficients of $f(t)g(t)$
are ratios of integers at most $nk$.

\medskip

Let us see the bit complexity of our algorithm in \S \ref{algorithm-sect}.
Fix $p$, $a,b,c$ and $\alpha$.
We need to compute the power series
\begin{equation}\label{bit-eq5}
\tau_{ab}^{(\sigma)}(t),\quad A^{(k)}_\sigma(t),\,\ldots,\, D^{(k)}_\sigma(t),
\quad E^{(k)}_\sigma(t)
\end{equation}
until the degree $pe_n+2p\sim p(p-1)/(p-2)n$.
First of all, the bit complexities of computing the constants
\[
\gamma_p+\psi_p(a^{(k)}), \quad\gamma_p+\psi_p(b^{(k)}), \quad
\log c
\]
modulo $p^n$ are small (cf. Appendix A), so that we can ignore them.
Moreover the power series $E^{(k)}_\sigma(t)$ is simple, so we can also ignore the bit complexity of computing it.

We observe the bit complexity of computing $\tau_{ab}^{(\sigma)}(t)$.
We work in a ring
\[
K[t]/(t^{pe_n+2p+1}),\quad K:=\Frac W.
\] 
We begin with the truncated polynomials
\[
F_{ab}(t)\in K[t]/(t^{pe_n+2p+1}).
\]
By \eqref{bit-eq1}, the bit complexities of computing all the coefficients are at most
\[
\sum_{i=0}^{pe_n+2p}O(i^2(\log i)^2)<O(n^3(\log n)^2).
\]
Next we need compute
\begin{equation}\label{bit-eq3}
\frac{1}{F_{ab}(t)}=(1+f)(1+f^2)\cdots(1+f^{2^d})
\in K[t]/(t^{pe_n+2p+1}),\quad f:=1-F_{ab}(t)
\end{equation}
where $d:=\lfloor\log_2(pe_n+2p)\rfloor+1\sim \log_2 n$.
The denominators and numerators of the coefficients of $f^k$ for $k\leq pe_n+2p$ 
are at most
\begin{equation}\label{bit-eq4}
\sum_{i_1+\cdots+i_k=l,\,i_r\geq 1}(i_1!\cdots i_k!)^2
<(l!)^2\binom{l-1}{k-1}<(l!)^2l^{pe_n+2p}
<(n!)^2n^{cn}
\end{equation}
with $c>0$ a constant.
Hence the bit complexities of computing $f^2,\ldots,f^{2^d}$ are at most
\[
O(n^2(\log (n!^2n^{cn}))^2)=O(n^4(\log n)^2)
\]
by \eqref{bit-eq2},
and hence the bit complexity of computing \eqref{bit-eq3} is
\[
O(dn^4(\log n)^2))=O(n^4(\log n)^3).
\]
Summing up the above, the bit complexity of computing $\tau^{(\sigma)}_{ab}(t)$
is $O(n^4(\log n)^3)$.

The power series of $A^{(k)}_\sigma(t),\,\ldots,\, D^{(k)}_\sigma(t)$
are obtained by applying standard arithmetic operations
(addition, subtraction and multiplication) 
on polynomials whose coefficients are ratios of integers at most \eqref{bit-eq4}.
Therefore the bit complexities do not exceed $O(n^4(\log n)^3)$.
All the algorithms in {\bf Step 1},\ldots, {\bf Step 3} are
standard arithmetic operations on the coefficients in the polynomials \eqref{bit-eq5}. 
One concludes that
the total bit complexity of the algorithm in \S \ref{algorithm-sect} is
$O(n^4(\log n)^3)$.


\section{Appendix A : $p$-adic polygamma functions}\label{AppendixA-sect}
We give a brief review of $p$-adic polygamma functions introduced in \cite[\S 2]{New}.

\medskip

Let $r\in\Z$ be an integer. For $z\in \Z_p$, define
\begin{equation}\label{wt-polygamma-def}
\wt{\psi}_p^{(r)}(z):=\lim_{n\in\Z_{>0},n\to z}\sum_{1\leq k<n,p\not{\hspace{0.7mm}|}\, k}
\frac{1}{k^{r+1}}
\end{equation}
where ``$n\to z$" means the limit with respect to the $p$-adic metric.
The existence of the limit follows from the fact that
\begin{equation}\label{equiv}
\sum_{1\leq k<p^n,p\not{\hspace{0.7mm}|}\,k}k^m\equiv 
\begin{cases}
-p^{n-1}&p\geq 3\mbox{ and }(p-1)|m\\
2^{n-1}&p=2\mbox{ and }2 |m\\
1&p=2\mbox{ and }n=1\\
0&\mbox{otherwise}
\end{cases}
\end{equation}
modulo $p^n$.
Thus $\wt\psi^{(r)}_p(z)$ is a $p$-adic continuous function on $\Z_p$.
Let $\log(x)$ be the Iwasawa logarithmic function which is characterized as the unique
continuous homomorphism $\log:\C^\times_p\to\C_p$ such that
$\log(x)=0$ if $x=p$ or $x$ is a root of unity, and 
\[
\log(x)=-\sum_{n=1}^\infty\frac{(1-x)^n}{n},\quad \text{if }|1-x|_p<1.
\]
Define the {\it $p$-adic Euler constant} by
\[
\gamma_p:=-\lim_{n\to\infty}\frac{1}{p^n}\sum_{0\leq j<p^n,p\not{\hspace{0.7mm}|}\,j}\log(j).\]
We define the $r$-th {\it $p$-adic polygamma function} to be
\begin{equation}\label{polygamma-def}
\psi_p^{(r)}(z):=\begin{cases}
-\gamma_p+\wt{\psi}^{(0)}_p(z)&r=0\\
-\zeta_p(r+1)+\wt{\psi}^{(r)}_p(z)&r\ne0
\end{cases}
\end{equation}
where $\zeta_p(r+1)$ is the special value of the $p$-adic zeta function
(see \cite[Lem 2.3]{New}).
If $r=0$, we also write $\psi_p(z)=\psi_p^{(0)}(z)$ and call it the {\it $p$-adic digamma function}.

\medskip

Concerning Dwork's $p$-adic hypergeometric functions,
we need to compute the special values of $\wt\psi_p(z)=\psi_p(z)+\gamma_p$
modulo $p^n$
(cf. \eqref{tau}).
To do this, the sum \eqref{wt-polygamma-def} is not useful because the number
of terms increases with exponential order by \eqref{equiv}.
However we can avoid this difficulty by using the following theorem.

\begin{thm}[{\cite[Thm.\,2.5]{New}}]\label{polygamma-thm2}
Let $0\leq i<N$ be integers and suppose $p\not|N$. Then
\begin{equation}\label{polygamma-thm2-eq}
\wt{\psi}_p^{(r)}\left(\frac{i}{N}\right)=N^r
\sum_{\ve\in \mu_N\setminus\{1\}}(1-\ve^{-i})\ln_{r+1}^{(p)}(\ve)
\end{equation}
where $\ln_k^{(p)}(z)$ are the $p$-adic polylogarithmic functions 
(cf. \cite[IV]{C-dlog}).
\end{thm}
Let $r=0$.
Then 
\begin{align*}
\ln_1^{(p)}(z)&=-p^{-1}\log\frac{(1-z)^p}{1-z^p}\\
&=\sum_{n=1}^\infty \frac{p^{n-1}}{n}w(z)^n,\quad
w(z):=p^{-1}\left(1-\frac{(1-z)^p}{1-z^p}\right).
\end{align*}
Using this expansion, one can compute 
$\wt\psi_p(i/N)$ mod $p^n$ without \eqref{wt-polygamma-def}.

\section{Appendix B : Resolution of Singularities}\label{AppendixB-sect}
Let $W$ be a commutative ring.
Let $X$ be a smooth $W$-scheme of relative dimension $d\geq 2$
or its completion along a closed subscheme.
A divisor $D$ is called a {\it relative normal crossing divisor} (abbreviated 
relative NCD)
over $W$
if it is locally
defined  by $x_1\cdots x_s=0$ where $(x_1,\ldots, x_d)$ is a local coordinates
over $W$.
Further $D$ is called {\it simple} if each component is smooth over $W$.

\begin{prop}\label{Appendix-1}
Let $n>0$ be an integer which is invertible in $W$.
Let
\begin{equation}\label{appA-eq1}
X:=\Spec W[[x,y,s]]/(sx-y^n)\supset C:=\Spec W[[x,y,s]]/(s,y).
\end{equation}
Then there is a proper morphism $\rho:X'\to X$ 
satisfying the following. Put $D:=\rho^{-1}(C)$.
\begin{itemize}
\item
$X'$ is smooth over $W$, and $X'\setminus D\os{\cong}{\to}
X\setminus C$,
\item
$D=
E_1+2E_2+\cdots+(n-1)E_{n-1}+nC'$ where $E_i$ are exceptional curves
and $C'$ is the proper transform of $C$,
\item
$E_1+E_2+\cdots+E_{n-1}+C'$ is a simple relative NCD over $W$. The figure is as follows.
\end{itemize}
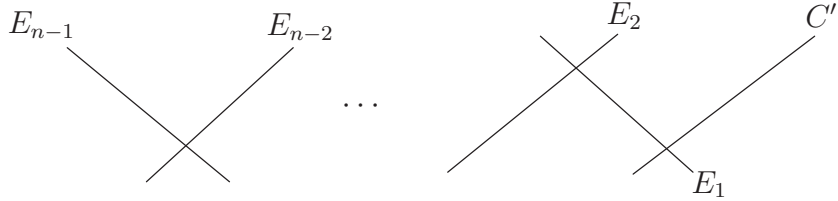
\begin{center}
{\unitlength 0.1in%
\begin{picture}(41.8000,9.4000)(1.7000,-11.0000)%
%
\special{pn 8}%
\special{pa 2450 1050}%
\special{pa 3330 330}%
\special{fp}%
\special{pa 2930 340}%
\special{pa 3720 1050}%
\special{fp}%
\special{pa 3410 1060}%
\special{pa 4350 340}%
\special{fp}%
%
\special{pn 8}%
\special{pa 480 400}%
\special{pa 1320 1100}%
\special{fp}%
\special{pa 890 1100}%
\special{pa 1650 400}%
\special{fp}%
\put(19.0000,-7.4000){\makebox(0,0)[lb]{$\cdots$}}%
\put(43.0000,-2.9000){\makebox(0,0)[lb]{$C'$}}%
\put(37.1000,-11.8000){\makebox(0,0)[lb]{$E_1$}}%
\put(15.1000,-3.7000){\makebox(0,0)[lb]{$E_{n-2}$}}%
\put(1.7000,-3.6000){\makebox(0,0)[lb]{$E_{n-1}$}}%
\put(32.7000,-3.0000){\makebox(0,0)[lb]{$E_2$}}%
\end{picture}}%
\end{center}
\end{prop}
\begin{pf}
Let $\rho_1:X_1\to X$ be the blow-up with center $(x,y,s)=(0,0,0)$.
Then $X_1$ is covered by affine open sets
\begin{align*}
U_1&=\Spec W[[x,y,s]][y_1,s_1]/(s_1-x^{n-2}y_1^n,xy_1-y,xs_1-s)\\
&\cong\Spec W[[x,y,s]][y_1]/(xy_1-y,x^{n-1}y_1^n-s),
\end{align*}
\begin{align*}
U_2&=\Spec W[[x,y,s]][x_2,y_2]/(x_2-s^{n-2}y_2^n,sy_2-y,sx_2-x)\\
&\cong\Spec W[[x,y,s]][y_2]/(sy_2-y,s^{n-1}y^n_2-x),
\end{align*}
\[
U_3=\Spec W[[x,y,s]][x_3,s_3]/(s_3x_3-y^{n-2},yx_3-x,ys_3-s).
\]
$U_1$ and $U_2$ are smooth over $W$.
If $n=2$, there is a unique exceptional curve $E$ such that $E\cap U_1=\{x=0\}$,
and $\rho_1^{-1}(C)=E+2C'$ where $C'$ is the proper transform
of $C$.
$X_1$ is smooth over $W$ and $E+C'$ is a simple relative NCD, so we are done.
If $n\geq 3$, then 
the divisor $D_1:=\rho_1^{-1}(C)=E_1+(n-1)E_2+nC'$ is as follows.

\begin{center}
{\unitlength 0.1in%
\begin{picture}( 31.4000, 12.7800)(  8.1000,-18.6800)%
\put(22.0000,-20.2000){\makebox(0,0)[lb]{$(n-1)E_2$}}%
\put(32.6000,-8.7000){\makebox(0,0)[lb]{$nC'$}}%
\put(8.1000,-20.2000){\makebox(0,0)[lb]{$E_1$}}%
\put(39.5000,-13.9000){\makebox(0,0)[lb]{Figure. $n\geq 3$}}%
%
\special{pn 8}%
\special{pa 1000 620}%
\special{pa 2350 1860}%
\special{fp}%
%
\special{pn 8}%
\special{pa 1020 1830}%
\special{pa 2370 590}%
\special{fp}%
%
\special{pn 8}%
\special{pa 1880 810}%
\special{pa 3190 800}%
\special{fp}%
%
\special{pn 4}%
\special{sh 1}%
\special{ar 1670 1232 16 16 0  6.28318530717959E+0000}%
\put(18.1000,-12.8000){\makebox(0,0)[lb]{$O$}}%
%
\special{pn 8}%
\special{ar 1672 1236 24 24  1.0040671  0.9600704}%
%
\special{pn 20}%
\special{pa 1674 1216}%
\special{pa 1674 1256}%
\special{fp}%
\special{pa 1686 1222}%
\special{pa 1686 1250}%
\special{fp}%
\special{pa 1662 1218}%
\special{pa 1662 1254}%
\special{fp}%
%
\special{pn 8}%
\special{ar 1673 1237 47 47  5.7288778  5.6084444}%
%
\special{pn 20}%
\special{pa 1667 1195}%
\special{pa 1667 1279}%
\special{fp}%
\special{pa 1679 1195}%
\special{pa 1679 1279}%
\special{fp}%
\special{pa 1691 1199}%
\special{pa 1691 1275}%
\special{fp}%
\special{pa 1703 1207}%
\special{pa 1703 1267}%
\special{fp}%
\special{pa 1655 1199}%
\special{pa 1655 1275}%
\special{fp}%
\special{pa 1643 1207}%
\special{pa 1643 1267}%
\special{fp}%
\end{picture}}%
\end{center}

Here $E_1$ and $E_2$ are exceptional curves such that
$E_1\cap U_3=\{y=x_3=0\}$ and $E_2\cap U_3=\{y=s_3=0\}$, and $O$ is the point
$(x_3,y,s_3)=(0,0,0)$ in $U_3$.
In a neighborhood of $O$, $X_1$ is locally defined by an equation $s_3x_3=y^{n-2}$.
If $n=3$, then $X_1$ is smooth, so we are done.
If $n\geq 4$, let $\rho_2:X_2\to X_1$ be the blow-up at $O$.

\begin{center}
{\unitlength 0.1in%
\begin{picture}(32.2000,11.0000)(5.6000,-16.3000)%
\put(31.4000,-8.0000){\makebox(0,0)[lb]{$4C'$}}%
\put(19.2000,-17.6000){\makebox(0,0)[lb]{$E_1$}}%
\put(6.0000,-17.6000){\makebox(0,0)[lb]{$3E_2$}}%
\put(22.6000,-13.6000){\makebox(0,0)[lb]{$2E_3$}}%
\put(37.8000,-12.9000){\makebox(0,0)[lb]{Figure. $n=4$}}%
%
\special{pn 8}%
\special{pa 1840 760}%
\special{pa 3080 760}%
\special{fp}%
%
\special{pn 8}%
\special{pa 2050 530}%
\special{pa 2050 1590}%
\special{fp}%
%
\special{pn 8}%
\special{pa 560 1310}%
\special{pa 2240 1310}%
\special{fp}%
%
\special{pn 8}%
\special{pa 760 530}%
\special{pa 760 1590}%
\special{fp}%
\end{picture}}%
\end{center}

If $n=4$, then $X_2$ is smooth over $W$ and $\rho_2^{-1}(D_1)=E_1+3E_2+2E_3+4C'$
is as in the figure where $E_3$ is the unique exceptional curve.
So we are done.
If $n\geq 5$, then $D_2=\rho_2^{-1}(D_1)=E_1+(n-1)E_2+2E_3+(n-2)E_4+nC'$ is as follows.

\begin{center}
{\unitlength 0.1in%
\begin{picture}( 31.1000, 10.6800)(  4.8000,-17.5800)%
\put(35.9000,-13.0000){\makebox(0,0)[lb]{Figure. $n\geq 5$}}%
\put(31.9000,-9.9000){\makebox(0,0)[lb]{$nC'$}}%
\put(24.1000,-19.1000){\makebox(0,0)[lb]{$E_1$}}%
\put(4.8000,-8.7500){\makebox(0,0)[lb]{$(n-1)E_2$}}%
\put(27.8000,-16.0000){\makebox(0,0)[lb]{$2E_3$}}%
\put(12.5000,-19.1000){\makebox(0,0)[lb]{$(n-2)E_4$}}%
\put(16.4000,-14.6600){\makebox(0,0)[lb]{$O$}}%
%
\special{pn 8}%
\special{pa 2480 690}%
\special{pa 2480 1750}%
\special{fp}%
%
\special{pn 8}%
\special{pa 1550 690}%
\special{pa 1550 1750}%
\special{fp}%
%
\special{pn 8}%
\special{pa 820 930}%
\special{pa 1840 930}%
\special{fp}%
%
\special{pn 8}%
\special{pa 2210 930}%
\special{pa 3130 930}%
\special{fp}%
%
\special{pn 4}%
\special{sh 1}%
\special{ar 1552 1532 16 16 0  6.28318530717959E+0000}%
%
\special{pn 8}%
\special{pa 1330 1535}%
\special{pa 2735 1540}%
\special{fp}%
%
\special{pn 8}%
\special{ar 1554 1534 44 44  0.3217506  0.2449787}%
%
\special{pn 20}%
\special{pa 1540 1496}%
\special{pa 1540 1572}%
\special{fp}%
\special{pa 1552 1494}%
\special{pa 1552 1574}%
\special{fp}%
\special{pa 1564 1494}%
\special{pa 1564 1574}%
\special{fp}%
\special{pa 1576 1500}%
\special{pa 1576 1568}%
\special{fp}%
\special{pa 1588 1514}%
\special{pa 1588 1554}%
\special{fp}%
\special{pa 1528 1504}%
\special{pa 1528 1564}%
\special{fp}%
\special{pa 1516 1524}%
\special{pa 1516 1546}%
\special{fp}%
\end{picture}}%
\end{center}

If $n=5$ then $X_2$ is smooth over $W$, and so we are done.
If $n\geq 6$, there is a singular point $O$.
In a neighborhood of $O$, $X_2$ is defined by an equation $s_4x_4=y^{n-4}$,
and $E_3=\{y=x_4=0\}$, $E_4=\{y=s_4=0\}$ and $D_2=\{s_4y^2=0\}$.
Then we take the blowing-up at $O$.
Continuing this, we finally obtain $\rho:X'=X_n\to X$  with $X'$ a smooth $W$-scheme 
such that $\rho^{-1}(C)=E_1+2E_2+\cdots+(n-1)E_{n-1}+nC'$ and
$E_1+\cdots+E_{n-1}+C'$ is a simple relative NCD over $W$.
\end{pf}
\begin{prop}\label{Appendix-2}
Let $N,M>0$ be integers which are invertible in $W$. 
Let $d=\gcd(N,M)$.
Suppose that $W$ contains a primitive $d$-th root of unity. Let
\begin{equation}\label{appA-eq2}
X:=\Spec W[[x,y]]\supset C:=\Spec W[[x,y]]/(x^N+y^M).
\end{equation}
Then there is a proper morphism $\rho:X'\to X$ satisfying the following.
Put $D:=\rho^{-1}(C)$.
\begin{itemize}
\item
$X'$ is smooth over $W$, and $X'\setminus D\os{\cong}{\to}
X\setminus C$,
\item
$D=\sum n_iD_i$ with $D_i$ smooth over $W$. Moreover 
$D=\sum D_i$ is a simple relative NCD over $W$, and
the multiplicities $n_i$
are either of
\[
1,\quad iN,\quad jM,\quad i\in\{1,\ldots,M\},\, j\in\{1,\ldots,N\}.
\]
\end{itemize}
The figure of $\sum_iD_i$ is as follows, where $C'$ is the proper transform of $C$ which
has $d$-components.
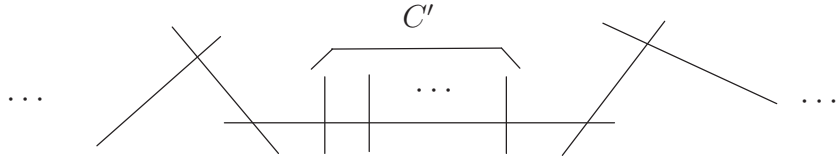
\begin{center}
{\unitlength 0.1in%
\begin{picture}(41.1000,8.0000)(2.0000,-11.5000)%
%
\special{pn 8}%
\special{pa 1330 980}%
\special{pa 3350 980}%
\special{fp}%
%
\special{pn 8}%
\special{pa 1850 740}%
\special{pa 1850 1140}%
\special{fp}%
%
\special{pn 8}%
\special{pa 2080 730}%
\special{pa 2080 1130}%
\special{fp}%
%
\special{pn 8}%
\special{pa 2790 740}%
\special{pa 2790 1140}%
\special{fp}%
%
\special{pn 8}%
\special{pa 1780 700}%
\special{pa 1890 595}%
\special{fp}%
\special{pa 1885 595}%
\special{pa 2765 590}%
\special{fp}%
\special{pa 2765 590}%
\special{pa 2860 695}%
\special{fp}%
\put(22.5000,-4.8000){\makebox(0,0)[lb]{$C'$}}%
\put(23.1500,-8.5000){\makebox(0,0)[lb]{$\cdots$}}%
%
\special{pn 8}%
\special{pa 3080 1150}%
\special{pa 3610 450}%
\special{fp}%
\special{pa 3290 460}%
\special{pa 4190 880}%
\special{fp}%
\put(43.1000,-9.1000){\makebox(0,0)[lb]{$\cdots$}}%
%
\special{pn 8}%
\special{pa 1610 1140}%
\special{pa 1060 480}%
\special{fp}%
\special{pa 1310 530}%
\special{pa 670 1100}%
\special{fp}%
\put(2.0000,-9.0000){\makebox(0,0)[lb]{$\cdots$}}%
\end{picture}}%
\end{center}

\end{prop}
\begin{pf}
For integers $a,b,c,d\geq 0$,
we denote by $I(a,b;c,d)$
the divisor $\Spec W[[x,y]]/(x^ay^b(x^c+y^d))$
in $\Spec W[[x,y]$.
Our goal is to compute the embedded resolution of $I(0,0;N,M)$.

Let $D=\Spec W[[x,y]]/(x^ay^b(x^c+y^d))=aD_x+bD_y+D_s\subset X$
where $D_x:=\{x=0\}$, $D_y:=\{y=0\}$ and $D_s:=\{x^c+y^d=0\}$.
Let $\rho:X'\to X$ be the blow-up with center $(x,y)=(0,0)$.
Then $D':=\rho^{-1}(D)=(a+b+c)E+aD_x'+bD_y'+D'_s$ where $D'_*$ denotes the proper
transform of $D_*$ and $E$ the exceptional curve. 
In case $c<d$, there is a unique point $O$ which is not normal crossing, and 
it is locally given by an equation $x^ay^{a+b+c}(x^c+y^{d-c})=0$, 
namely $I(a,a+b+c;c,d-c)$.
The multiplicities of $D'$ are $1,a,b,a+b+c$.
In case $c>d$, there is also a unique point $O$ such that
the divisor $D'$ around $O$ is $I(a+b+d,b;c-d,d)$, and
the multiplicities of $D_1$ are $1,a,b,a+b+d$.
In case $c=d$, the divisor $D'=(a+b+c)E+aD_x'+bD_y'+D'_s$
satisfies that
$E+D_x'+D_y'+D'_s$ is a simple relative NCD over $W$, and
$D'_s$ has $c$-components (see the figure).
In this case we stop the resolution.
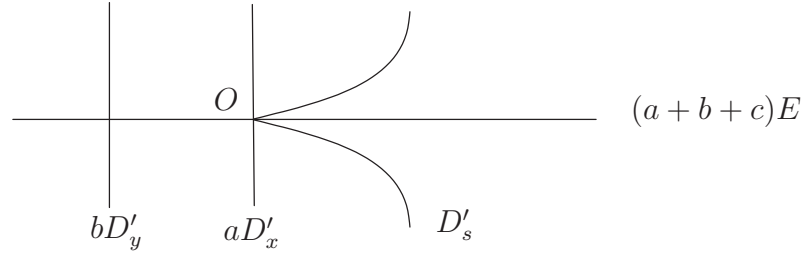
\begin{center}
{\unitlength 0.1in%
\begin{picture}(33.4000,14.4400)(1.1000,-18.7200)%
%
\special{pn 8}%
\special{pa 250 1310}%
\special{pa 3270 1310}%
\special{fp}%
%
\special{pn 8}%
\special{pa 750 700}%
\special{pa 750 1770}%
\special{fp}%
%
\special{pn 8}%
\special{pa 1490 710}%
\special{pa 1500 1760}%
\special{fp}%
\put(34.5000,-13.5000){\makebox(0,0)[lb]{$(a+b+c)E$}}%
\put(13.4000,-19.8000){\makebox(0,0)[lb]{$aD'_x$}}%
\put(6.5000,-19.9000){\makebox(0,0)[lb]{$bD'_y$}}%
\put(1.1000,-5.8000){\makebox(0,0)[lb]{Case $c<d$}}%
%
\special{pn 8}%
\special{pa 1496 1308}%
\special{pa 1534 1299}%
\special{pa 1610 1279}%
\special{pa 1648 1270}%
\special{pa 1686 1260}%
\special{pa 1723 1250}%
\special{pa 1831 1217}%
\special{pa 1866 1205}%
\special{pa 1900 1193}%
\special{pa 1933 1181}%
\special{pa 1997 1153}%
\special{pa 2027 1139}%
\special{pa 2056 1123}%
\special{pa 2084 1107}%
\special{pa 2110 1090}%
\special{pa 2135 1072}%
\special{pa 2159 1053}%
\special{pa 2181 1033}%
\special{pa 2202 1012}%
\special{pa 2221 990}%
\special{pa 2238 967}%
\special{pa 2253 943}%
\special{pa 2267 917}%
\special{pa 2278 891}%
\special{pa 2288 862}%
\special{pa 2295 833}%
\special{pa 2300 802}%
\special{pa 2304 770}%
\special{pa 2306 748}%
\special{fp}%
%
\special{pn 8}%
\special{pa 1496 1312}%
\special{pa 1534 1321}%
\special{pa 1610 1341}%
\special{pa 1648 1350}%
\special{pa 1686 1360}%
\special{pa 1723 1370}%
\special{pa 1831 1403}%
\special{pa 1866 1415}%
\special{pa 1900 1427}%
\special{pa 1933 1439}%
\special{pa 1997 1467}%
\special{pa 2027 1481}%
\special{pa 2056 1497}%
\special{pa 2084 1513}%
\special{pa 2110 1530}%
\special{pa 2135 1548}%
\special{pa 2159 1567}%
\special{pa 2181 1587}%
\special{pa 2202 1608}%
\special{pa 2221 1630}%
\special{pa 2238 1653}%
\special{pa 2253 1677}%
\special{pa 2267 1703}%
\special{pa 2278 1729}%
\special{pa 2288 1758}%
\special{pa 2295 1787}%
\special{pa 2300 1818}%
\special{pa 2304 1850}%
\special{pa 2306 1872}%
\special{fp}%
\put(24.4000,-19.4000){\makebox(0,0)[lb]{$D'_s$}}%
\put(12.9000,-12.6000){\makebox(0,0)[lb]{$O$}}%
\end{picture}}%
\end{center}
\vspace{1cm}
\begin{center}
{\unitlength 0.1in%
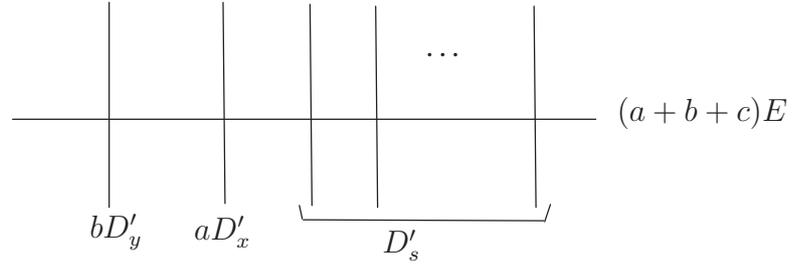
\begin{picture}(32.8000,15.4000)(1.0000,-19.0500)%
%
\special{pn 8}%
\special{pa 250 1310}%
\special{pa 3270 1310}%
\special{fp}%
%
\special{pn 8}%
\special{pa 750 700}%
\special{pa 750 1770}%
\special{fp}%
%
\special{pn 8}%
\special{pa 1340 700}%
\special{pa 1350 1750}%
\special{fp}%
\put(33.8000,-13.6000){\makebox(0,0)[lb]{$(a+b+c)E$}}%
\put(11.9000,-19.8000){\makebox(0,0)[lb]{$aD'_x$}}%
\put(6.5000,-19.9000){\makebox(0,0)[lb]{$bD'_y$}}%
\put(1.0000,-5.1000){\makebox(0,0)[lb]{Case $c=d$}}%
%
\special{pn 8}%
\special{pa 1790 710}%
\special{pa 1800 1760}%
\special{fp}%
%
\special{pn 8}%
\special{pa 2130 720}%
\special{pa 2140 1770}%
\special{fp}%
%
\special{pn 8}%
\special{pa 2950 720}%
\special{pa 2960 1770}%
\special{fp}%
\put(23.9000,-9.8000){\makebox(0,0)[lb]{\ldots}}%
\put(21.6500,-20.5000){\makebox(0,0)[lb]{$D'_s$}}%
%
\special{pn 8}%
\special{pa 1755 1835}%
\special{pa 3000 1840}%
\special{fp}%
%
\special{pn 8}%
\special{pa 1750 1830}%
\special{pa 1735 1760}%
\special{fp}%
%
\special{pn 8}%
\special{pa 3005 1845}%
\special{pa 3035 1755}%
\special{fp}%
\end{picture}}%
\end{center}
Define
\begin{equation}\label{Euclid-eq1}
(I(a,b;c,d))':=\begin{cases}
I(a,a+b+c;c,d-c)&c\leq d\\
I(a+b+d,b;c-d,d)&c>d\\
I(a,b;c,d)&cd=0
\end{cases}
\end{equation}
and $I^{(0)}=I$, $I^{(i)}=(I^{(i-1)})'$.
We begin with $I(0,0;N,M)$ and 
consider a sequence $I(a_i,b_i;c_i,d_i):=(I(0,0;N,M))^{(i)}$
\[
I(0,0;N,M),\,
I(a_1,b_1;c_1,d_1),\,\ldots,\,
I(a_{n},b_{n};c_{n},d_{n})
\]
until $c_nd_n=0$.
This corresponds to the sequence of blowing ups at $O$'s as above
\[
X_n\lra X_{n-1}\lra \cdots\lra X_0=X
\]
such that the inverse image of $C$ in $X_n$ is supported in a relative simple NCD.
Moreover let $D_i\subset X_i$ be the inverse image of $C$.
Then the multiplicities of $D_i$ are either of $1,a_1,\ldots,a_i,b_1,\ldots,b_i$.
Therefore if we show Lemma \ref{Appendix-lem} below (which is a simple lemma in
elementary number theory), then it ends the proof of Proposition \ref{Appendix-2}.
\end{pf}

\begin{lem}\label{Appendix-lem}
Let $N,M\geq 1$ be integers and let 
$I(a_i,b_i;c_i,d_i):=(I(0,0,;N,M))^{(i)}$ be defined by \eqref{Euclid-eq1}.
Let $n$ be the minimal integer such that $c_nd_n=0$.
\begin{enumerate}
\item[$(1)$]
There are integers $A_i,B_i,C_i,D_i\geq0$ such that
$a_i=A_iM$, $b_i=B_iN$, $c_i=C_iN-A_iM$, $d_i=D_iM-B_iN$.
\item[$(2)$]
$A_i$, $B_i$, $C_i$, $D_i$ are non-decreasing sequences, and 
$A_n,D_n\leq N$ and $B_n,C_n\leq M$.
\end{enumerate}
\end{lem}
\begin{pf}
(1)
The assertion is clear for $i=0$
by putting $(A_0,B_0,C_0,D_0):=(0,0,1,1)$.
Suppose that the assertion holds for $i$.
By definition
\begin{align*}
(a_{i+1},b_{i+1},c_{i+1},d_{i+1})
&=
\begin{cases}
(a_i,a_i+b_i+c_i,c_i,d_i-c_i)&c_i\leq d_i\\
(a_i+b_i+d_i,b_i,c_i-d_i,d_i)&c_i>d_i
\end{cases}\\
&=
\begin{cases}
(A_iM,(B_i+C_i)N,C_iN-A_iM,(A_i+D_i)M-(B_i+C_i)M)\\
((A_i+D_i)M,B_iN,(B_i+C_i)N-(A_i+D_i)M,D_iM-B_iN).
\end{cases}
\end{align*}
Hence the assertion holds by putting
\begin{equation}\label{Appendix-lem-eq1}
(A_{i+1},B_{i+1},C_{i+1},D_{i+1})
:=
\begin{cases}
(A_i,B_i+C_i,C_i,A_i+D_i)&c_i\leq d_i\\
(A_i+D_i,B_i,B_i+C_i,D_i)&c_i>d_i.
\end{cases}
\end{equation}

\noindent
(2) The former assertion is obvious from \eqref{Appendix-lem-eq1}.
We show $A_n,D_n\leq N$ and $B_n,C_n\leq M$.
The algorithm $(c_0,d_0)\to (c_1,d_1)\to \cdots\to (c_n,d_n)$ is the Euclidean algorithm.
Therefore $(c_n,d_n)=(0,\gcd(N,M))$ or $(\gcd(N,M),0)$.
In case $(c_n,d_n)=(0,\gcd(N,M))$, 
$A_n,B_n,C_n,D_n$ are characterized as the minimal positive 
integers satisfying
$C_nN=A_nM$ and $D_nM-B_nN=\gcd(N,M)$. 
Hence it turns out that $A_n,D_n\leq N$ and $B_n,C_n\leq M$.
The conclusion is the same also in case $(c_n,d_n)=(\gcd(N,M),0)$.
\end{pf}


\noindent
Department of Mathematics, Hokkaido University,
\par\noindent
Sapporo 060-0810,
JAPAN

\smallskip

\noindent
{\it E-mail} : \textbf{asakura@math.sci.hokudai.ac.jp}


\end{document}